\documentclass[final]{siamonline171218}

\usepackage{mathtools} 
\usepackage{thmtools} 
\usepackage{hyperref} 
\usepackage{color} 
\usepackage{amssymb} 
\usepackage{tikz} 
\usepgflibrary{shapes.misc}
\usetikzlibrary{positioning}
\usetikzlibrary{matrix} 
\usetikzlibrary{patterns}
\usepackage{mleftright} 
\usepackage{bbm} 
\usepackage{enumitem}
\usepackage{cleveref}

\crefname{convariable}{Condition}{Conditions}
\Crefname{convariable}{Condition}{Conditions}
\crefname{probvariable}{Problem}{Problems}
\Crefname{probvariable}{Problem}{Problems}
\crefname{myprobsec}{Problem}{Problems}
\Crefname{myprobsec}{Problem}{Problems}
\crefname{step}{Step}{Steps}
\Crefname{step}{Step}{Steps}
\crefname{point}{Point}{Points}
\Crefname{point}{Point}{Points}
\crefname{part}{Part}{Parts}
\Crefname{part}{Part}{Parts}
\crefname{condition}{Condition}{Conditions}
\Crefname{condition}{Condition}{Conditions}
\crefname{problem}{Problem}{Problems}
\Crefname{problem}{Problem}{Problems}

\setlist[enumerate]{leftmargin=.5in}
\setlist[itemize]{leftmargin=.5in}

\headers{The Helmholtz equation in random media}{O.~R.~Pembery and E.~A.~Spence}

\ifpdf
\hypersetup{
  pdftitle={The Helmholtz equation in random media: well-posedness and a priori bounds},
  pdfauthor={O. R. Pembery and E. A. Spence}
}
\fi

\newcommand{\Li}[1]{L^{\infty}\mleft(#1\mright)}

\newcommand{\LtD}{\Lt{D}}

\newcommand{\LtOLtD}{\Lt{\Omega;\Lt{D}}}
\newcommand{\Lt}[1]{L^{2}\mleft(#1\mright)}
\newcommand{\OFP}{\mleft(\Omega,\cF,\PP\mright)}

\newcommand{\TC}{\cT_{\cC}}

\newcommand{\abs}[1]{\mleft|#1\mright|}
\newcommand{\biggabs}[1]{\bigg|#1\bigg|}
\newcommand{\Bigabs}[1]{\Big|#1\Big|}

\newcommand{\bas}{\begin{assumption}}
\newcommand{\bde}{\begin{definition}}
\newcommand{\bprob}{\begin{problem}}
\newcommand{\bx}{\mathbf{x}}
\newcommand{\dd}{\,\mathrm{d}}
\newcommand{\eas}{\end{assumption}}
\newcommand{\ede}{\end{definition}}
\newcommand{\eprob}{\end{problem}}
\newcommand{\grad}{\nabla}

\newcommand{\HmhGR}{\Hmh{\GR}}
\newcommand{\Hmh}[1]{H^{-1/2}\mleft(#1\mright)}

\newcommand{\set}[1]{\mleft\{#1\mright\}}

\newcommand{\tforall}{\text{ for all }}

\newcommand{\Acomega}{\cA_{c\mleft(\omega\mright)}}
\newcommand{\Lcomega}{\rhs_{c\mleft(\omega\mright)}}
\newcommand{\LtOX}{\Lt{\Omega;X}}
\newcommand{\LtOY}{\Lt{\Omega;Y}}
\newcommand{\LtOYas}{\LtOY^{*}}
\newcommand{\sol}{\cS}
\newcommand{\coeff}{c}

\newcommand{\pair}{P}
	
\newcommand{\rhs}{\cL}

\newcommand{\toform}{\cA}

\newcommand{\homspace}{\mathrm{B}\mleft(X,\Ys\mright)}

\newcommand{\NYs}[1]{\N{#1}_{\Ys}}

\newcommand{\Ys}{Y^{*}}
\newcommand{\Nhomspace}[1]{\N{#1}_{\homspace}}
\newcommand{\NX}[1]{\N{#1}_{X}}

\newcommand{\NLtOX}[1]{\N{#1}_{\LtOX}}

\newcommand{\dinfty}{d_{\infty}}

\newcommand{\NLtD}[1]{\N{#1}_{\LtD}}

\newcommand{\SAc}{\SA}
\newcommand{\AAA}{\mathfrak{A}}
\newcommand{\SLc}{\SL}
\newcommand{\LLL}{\mathfrak{L}}
\newcommand{\ud}{u_{d}}
\newcommand{\vn}{v_{n}}
\newcommand{\LpOX}{\Lp{\Omega;X}}
\newcommand{\Lp}[1]{L^{p}\mleft(#1\mright)}

\newcommand{\LoO}{\Lo{\Omega}}
\newcommand{\Lo}[1]{L^{1}\mleft(#1\mright)}

\newcommand{\NLoO}[1]{\N{#1}_{\LoO}}

\newcommand{\defn}[1]{\emph{#1}}
\newcommand{\fullrhs}{\rhs\circ\coeff}
\newcommand{\NY}[1]{\N{#1}_{Y}}
\newcommand{\Borel}[1]{\cB\mleft(#1\mright)}
\newcommand{\bigBorel}[1]{\cB\big(#1\big)}

\newcommand{\Ball}[3]{B^{#1}_{#2}\mleft(#3\mright)}

\newcommand{\LtOYs}{L^{2}\mleft(\Omega;Y^{*}\mright)}
\newcommand{\NLtO}[1]{\N{#1}_{\LtO}}

\newcommand{\LtO}{\Lt{\Omega}}
\newcommand{\NLtOYs}[1]{\N{#1}_{\LtOYs}}
\newcommand{\NLtOY}[1]{\N{#1}_{\LtOY}}
\newcommand{\SA}{\AAA}
\newcommand{\SL}{\LLL}
\newcommand{\LiOhomspace}{\Li{\Omega;\homspace}}
\newcommand{\vo}{v_{1}}
\newcommand{\vt}{v_{2}}
\newcommand{\half}{\frac{1}{2}}

\DeclareMathOperator{\esssup}{ess\,sup}
\DeclareMathOperator{\essinf}{ess\,inf}
\newcommand{\NLiOhomspace}[1]{\N{#1}_{\LiOhomspace}}
\newcommand{\HoD}{\Ho{D}}
\newcommand{\Ho}[1]{H^{1}\mleft(#1\mright)}

\newcommand{\IP}[2]{\mleft\langle#1,#2\mright\rangle}

\newcommand{\IPbig}[2]{\big\langle#1,#2\big\rangle}

\newcommand{\RRd}{\RR^{d}}

\newcommand{\LtOHoD}{\Lt{\Omega;\HoD}}

\newcommand{\Lq}[1]{L^{q}\mleft(#1\mright)}
\DeclareMathOperator{\supp}{ess\,supp}

\newcommand{\bxi}{\boldsymbol{\xi}}

\newcommand{\tendi}{\rightarrow \infty}

\newcommand{\NLtOLtD}[1]{\N{#1}_{\LtOLtD}}

\newcommand{\N}[1]{\mleft\|#1\mright\|}

\newcommand{\Top}[1]{\cT_{#1}}
\newcommand{\eps}{\varepsilon}
\newcommand{\de}{\coloneqq}
\newcommand{\st}{:}

\newcommand{\DR}{D_{R}}
\newcommand{\Dp}{D_{+}}
\newcommand{\BR}{B_{R}}
\newcommand{\HozDDR}{H^{1}_{0,D}\mleft(\DR\mright)}
\newcommand{\LtDp}{\Lt{\Dp}}
\newcommand{\compcont}{\subset\subset}
\newcommand{\LiDpRR}{\Li{\Dp;\RR}}
\newcommand{\nmin}{n_{\min}}
\newcommand{\nmax}{n_{\max}}

\newcommand{\Amin}{A_{\min}}
\newcommand{\vtb}{\overline{\vt}}
\newcommand{\LtOHozDDR}{\Lt{\Omega;\HozDDR}}
\newcommand{\as}{\mathfrak{a}}

\newcommand{\Czo}[1]{C^{0,1}\mleft(#1\mright)}

\newcommand{\WoiDRRRdtd}{W^{1,\infty}\mleft(\DR;\RR^{d\times d}\mright)}
\newcommand{\WoiDRRRdtdsmall}{W^{1,\infty}(\DR;\RR^{d\times d})}
\newcommand{\LiDRRR}{\Li{\DR;\RR}}
\newcommand{\LtDR}{\Lt{\DR}}

\newcommand{\NLiDRRRdtd}[1]{\N{#1}_{\LiDRRRdtd}}
\newcommand{\LiDRRRdtd}{\Li{\DR;\RR^{d\times d}}}
\newcommand{\NLiDRRR}[1]{\N{#1}_{\LiDRRR}}
\newcommand{\NLtDR}[1]{\N{#1}_{\LtDR}}

\newcommand{\IPGR}[2]{\IP{#1}{#2}_{\GR}}

\newcommand{\IPGRbig}[2]{\IPbig{#1}{#2}_{\GR}}

\newcommand{\GR}{\Gamma_{R}}
\newcommand{\Nw}[1]{\N{#1}_{1,k}}
\newcommand{\LtOLtDR}{\Lt{\Omega;\LtDR}}
\newcommand{\LiOLiDRRRdtd}{\Li{\Omega;\LiDRRRdtd}}
\newcommand{\LiOLiDRRR}{\Li{\Omega;\LiDRRR}}

\newcommand{\Woi}[1]{W^{1,\infty}\mleft(#1\mright)}
\newcommand{\RRdtd}{\RR^{d\times d}}
\newcommand{\Ls}{\mathfrak{L}}
\newcommand{\NLtOHozDDR}[1]{\N{#1}_{\LtOHozDDR}}
\newcommand{\BorelOX}{\cB\mleft(\Omega,X\mright)}

\newcommand{\Ao}{A_{1}}
\newcommand{\no}{n_{1}}
\newcommand{\fo}{f_{1}}
\newcommand{\Am}{A_{m}}
\newcommand{\nm}{n_{m}}
\newcommand{\fm}{f_{m}}
\newcommand{\Az}{A_{0}}
\newcommand{\nz}{n_{0}}
\newcommand{\fz}{f_{0}}
\newcommand{\DtN}{T_{R}}

\newcommand{\NHhGR}[1]{\N{#1}_{\HhGR}}
\newcommand{\HhGR}{\Hh{\GR}}
\newcommand{\Hh}[1]{H^{1/2}\mleft(#1\mright)}
\newcommand{\uz}{u_{0}}
\newcommand{\uo}{u_{1}}

\newcommand{\HozDDRas}{\mleft(\HozDDR\mright)^{*}}
\newcommand{\NHozDDRas}[1]{\N{#1}_{\HozDDRas}}
\newcommand{\CzoDR}{\Czo{\DR}}
\newcommand{\WoiDR}{\Woi{\DR}}
\newcommand{\Dm}{D_{-}}
\newcommand{\GD}{\Gamma_{D}}

\newcommand{\muo}{\mu_{1}}
\newcommand{\mut}{\mu_{2}}
\newcommand{\Co}{C_{1}}
\newcommand{\TS}{\Top{S}}

\newcommand{\evalbase}{\eta}

\newcommand{\HoDR}{\Ho{\DR}}
\newcommand{\ut}{u_{2}}
\newcommand{\owenshift}{0.5}
\newcommand{\owentextshift}{0.125}

\newcommand{\M}{M}
\newcommand{\FM}{\cM}
\newcommand{\FN}{\cN}
\newcommand{\TopT}{\Top{T}}
\newcommand{\NV}[1]{\N{#1}_{V}}
\newcommand{\NLpOV}[1]{\N{#1}_{\LpOV}}
\newcommand{\LpOV}{\Lp{\Omega;V}}
\newcommand{\NLiOV}[1]{\N{#1}_{\LiOV}}
\newcommand{\LiOV}{\Li{\Omega;V}}
\newcommand{\D}{\cD}
\newcommand{\Ind}[1]{\mathbbm{1}_{#1}}

\newcommand{\diff}{\delta}
\newcommand{\Omegay}{\Omega_{y}}

\newcommand{\Omegat}{\widetilde{\Omega}}
\newcommand{\yn}{y_{n}}
\newcommand{\Omegayn}{\Omega_{\yn}}
\newcommand{\Omegaync}{\Omega^{c}_{\yn}}
\newcommand{\ynmseq}{\mleft(\ynm\mright)_{m \in \NN}}
\newcommand{\ynm}{y_{n_{m}}}
\newcommand{\Omegaynm}{\Omega_{\ynm}}

\newcommand{\LtORR}{\Lt{\Omega;\RR}}
\newcommand{\TrR}{T_{R}}

\newcommand{\Xo}{X_{1}}
\newcommand{\Xm}{X_{m}}
\newcommand{\dm}{d_{m}}
\newcommand{\done}{d_{1}}
\newcommand{\xo}{x_{1}}
\newcommand{\xm}{x_{m}}
\newcommand{\yo}{y_{1}}
\newcommand{\ym}{y_{m}}
\newcommand{\dmetj}{d_{j}}
\newcommand{\xj}{x_{j}}
\newcommand{\yj}{y_{j}}
\newcommand{\aGm}{a_{\Az,\nz}}
\newcommand{\Lh}{L_{\fz}}
\newcommand{\Prodf}{P}
\newcommand{\measrectmany}[1]{\textrm{Rect}\mleft(#1\mright)}
\newcommand{\cCo}{\cC_{1}}
\newcommand{\cCm}{\cC_{m}}
\newcommand{\cCj}{\cC_{j}}
\newcommand{\Ctilde}{\widetilde{C}}
\newcommand{\nzmin}{n_{0,\min}}

\newcommand{\Azmin}{A_{0,\min}}

\newcommand{\CCd}{\CC^{d}}

\newcommand{\consto}{\tau_{1}}
\newcommand{\constt}{\tau_{2}}

\newcommand{\Bo}{B_{1}}
\newcommand{\Bt}{B_{2}}
\newcommand{\ft}{f_{2}}

\newcommand{\pDm}{\GD}
\newcommand{\uI}{u_{I}}
\newcommand{\LfO}{\Lf{\Omega}}
\newcommand{\Lf}[1]{L^{4}\mleft(#1\mright)}
\newcommand{\LfOLtD}{\Lf{\Omega;\LtD}}
\newcommand{\ri}{\mathrm{i}} 
\newcommand{\NLtOLtDR}[1]{\N{#1}_{\LtOLtDR}}
\newcommand{\bxz}{\mathbf{\xz}}
\newcommand{\xz}{x_{0}}
\newcommand{\Eo}{E_{1}}
\newcommand{\Et}{E_{2}}
\newcommand{\Eth}{E_{3}}

\newcommand{\Yj}{Y_{j}}
\newcommand{\psij}{\psi_{j}}
\newcommand{\Unif}{\mathrm{Unif}}
\newcommand{\bzero}{\mathbf{0}}

\newcommand{\WoiDRRR}{\Woi{\DR;\RR}}
\newcommand{\Zj}{Z_{j}}
\newcommand{\Psij}{\Psi_{j}}
\newcommand{\NopCCd}[1]{\N{#1}_{2}}

\newcommand{\LiDpRRdtd}{\Li{\Dp;\RRdtd}}
\newcommand{\Amax}{A_{\max}}
\newcommand{\Cj}{C_{j}}
\newcommand{\fj}{f_{j}}

\newcommand{\clos}[1]{\overline{#1}}

\newcommand{\NW}[1]{\N{#1}_{1,k}}
\newcommand{\kz}{k_{0}}

\newcommand{\tiff}{\text{ iff }}

\newcommand{\NHoD}[1]{\N{#1}_{\HoD}}
\newcommand{\LqOYas}{\LqOY^{*}}
\newcommand{\LqOY}{\Lq{\Omega;Y}}

\newcommand{\scalarmult}{M_{v}}
\newcommand{\inc}{\iota}
\DeclareMathOperator{\spansetsub}{span}
\newcommand{\spanset}[1]{\spansetsub\set{#1}}
\newcommand{\QQFF}{\QQ_{\FF}}
\newcommand{\phio}{\phi_{1}}
\newcommand{\phit}{\phi_{2}}
\newcommand{\QQU}{\QQ_{U}}
\newcommand{\setsq}{S_{s,q}}
\newcommand{\NU}[1]{\N{#1}_{U}}
\newcommand{\dw}{r_{\omega}}
\newcommand{\Dmclos}{\clos{\Dm}}

\newcommand{\NWoiDRRRdtd}[1]{\N{#1}_{\WoiDRRRdtdsmall}}

\newcommand{\LiRDRRR}{L^{\infty}_{R}\mleft(\DR;\RR\mright)}
\newcommand{\LiRminDRRR}{L^{\infty}_{R,\min}\mleft(\DR;\RR\mright)}
\newcommand{\LtRDR}{L^{2}_{R}\mleft(\DR\mright)}
\newcommand{\Leb}{\lambda}
\newcommand{\Real}[1]{\mathfrak{R}\mleft(#1\mright)}

\newcommand{\NLfOLtDR}[1]{\N{#1}_{\LfOLtDR}}
\newcommand{\LfOLtDR}{\Lf{\Omega;\LtDR}}

\newcommand{\LiRDRRRdtd}{L^{\infty}_{R}\mleft(\DR;\RRdtd\mright)}
\newcommand{\LiRminDRRRdtd}{L^{\infty}_{R,\min}\mleft(\DR;\RRdtd\mright)}

\newcommand{\WoiRminDRRRdtd}{W^{1,\infty}_{R,\min}\mleft(\DR;\RRdtd\mright)}

\newcommand{\NTA}[1]{\mathrm{NT}_{A}\mleft(#1\mright)}
\newcommand{\NTn}[1]{\mathrm{NT}_{n}\mleft(#1\mright)}
\newcommand{\Do}{D_{1}}
\newcommand{\Dt}{D_{2}}
\newcommand{\LiDp}{\Li{\Dp}}
\newcommand{\LiDpCC}{\Li{\Dp;\CC}}

\newcommand{\bxib}{\overline{\bxi}}
\newcommand{\alphanz}{\alpha_{\nz}}
\newcommand{\alphaAz}{\alpha_{\Az}}
\newcommand{\NLtOHoD}[1]{\N{#1}_{\LtOHoD}}
\newcommand{\NWoiDRRR}[1]{\N{#1}_{\WoiDRRR}}
\newcommand{\metsp}{W}
\newcommand{\metric}{d}
\newcommand{\Omegahat}{\hat{\Omega}}
\newcommand{\jhat}{\hat{j}}
\newcommand{\Yjhat}{Y_{\jhat}}
\newcommand{\Omegatilde}{\widetilde{\Omega}}
\newcommand{\zetatilde}{\zeta}
\newcommand{\LtOHozDDRnormal}{L^2(\Omega,\HozDDR)}

\newcommand{\beq}{\begin{equation}}
\newcommand{\eeq}{\end{equation}}
\newcommand{\beqs}{\begin{equation*}}
\newcommand{\eeqs}{\end{equation*}}
\newcommand{\bit}{\begin{itemize}}
\newcommand{\eit}{\end{itemize}}
\newcommand{\ben}{\begin{enumerate}}
\newcommand{\een}{\end{enumerate}}
\newcommand{\bal}{\begin{align}}
\newcommand{\eal}{\end{align}}
\newcommand{\bals}{\begin{align*}}
\newcommand{\eals}{\end{align*}}
\newcommand{\bse}{\begin{subequations}}
\newcommand{\ese}{\end{subequations}}
\newcommand{\bpr}{\begin{proposition}}
\newcommand{\epr}{\end{proposition}}
\newcommand{\bre}{\begin{remark}}
\newcommand{\ere}{\end{remark}}
\newcommand{\bpf}{\begin{proof}}
\newcommand{\epf}{\end{proof}}
\newcommand{\ble}{\begin{lemma}}
\newcommand{\ele}{\end{lemma}}
\newcommand{\bco}{\begin{corollary}}
\newcommand{\eco}{\end{corollary}}
\newcommand{\bex}{\begin{example}}
\newcommand{\eex}{\end{example}}
\newcommand{\bth}{\begin{theorem}}
\newcommand{\enth}{\end{theorem}}
\newcommand{\bcon}{\begin{condition}}
\newcommand{\econ}{\end{condition}}

\newcommand{\ton}{\text{ on }}
\newcommand{\tin}{\text{ in }}
\newcommand{\tfa}{\text{ for all }}
\newcommand{\tfor}{\text{ for }}

\newcommand{\tand}{\text{ and }}

\newcommand{\tif}{\text{ if }}

\newcommand{\spd}{\mathsf{SPD}}

\newcommand{\Zo}{Z_{1}}
\newcommand{\Zt}{Z_{2}}
\newcommand{\homspaceZoZt}{\mathrm{B}\mleft(\Zo,\Zt\mright)}
\newcommand{\homspaceYF}{\mathrm{B}\mleft(Y,\FF\mright)}
\newcommand{\evalZoZt}{\evalbase_{\Zo,\Zt}}
\newcommand{\Pmap}{\cP}
\newcommand{\toprodPv}{\pair_{\Pmap,v}}
\newcommand{\compPv}{\pi_{\Pmap,v}}

\newcommand{\NLiDRspd}[1]{\N{#1}_{\LiDRspd}}
\newcommand{\LiDRspd}{\Li{\DR;\spd}}

\newcommand{\WoiDRspd}{W^{1,\infty}(\DR;\spd)}
\newcommand{\NWoiDRspd}[1]{\N{#1}_{\WoiDRspd}}

\newcommand{\CC}{\mathbb C}

\newcommand{\FF}{\mathbb F}

\newcommand{\NN}{\mathbb N}

\newcommand{\PP}{\mathbb P}
\newcommand{\QQ}{\mathbb Q}
\newcommand{\RR}{\mathbb R}

\newcommand{\cA}{\mathcal A}
\newcommand{\cB}{\mathcal B}
\newcommand{\cC}{\mathcal C}
\newcommand{\cD}{\mathcal D}

\newcommand{\cF}{\mathcal F}
\newcommand{\cG}{\mathcal G}
\newcommand{\cH}{\mathcal H}

\newcommand{\cL}{\mathcal L}
\newcommand{\cM}{\mathcal M}
\newcommand{\cN}{\mathcal N}

\newcommand{\cP}{\mathcal P}

\newcommand{\cS}{\mathcal S}
\newcommand{\cT}{\mathcal T}

\title{The Helmholtz equation in random media: well-posedness and a priori bounds}

\author{O.~R.~Pembery\thanks{Department of Mathematical Sciences, University of Bath, Claverton Down, Bath, BA2 7AY, UK. This author is supported by a scholarship from the EPSRC Centre for Doctoral Training in Statistical Applied Mathematics at Bath (SAMBa), under the project EP/L015684/1. (\email{O.R.Pembery@bath.ac.uk}).}
\and E.~A.~Spence\thanks{Department of Mathematical Sciences, University of Bath, Claverton Down, Bath, BA2 7AY, UK. This author is supported by EPSRC grant EP/R005591/1. (\email{E.A.Spence@bath.ac.uk}).}}

\begin{document}

\maketitle

\numberwithin{equation}{section}

\newsiamthm{assumption}{Assumption}
\newsiamthm{example}{Example}
\newsiamremark{remark}{Remark}
\newsiamremark{condition}{Condition}
\newsiamremark{problem}{Problem}
\newsiamremark{neverused}{Labelled Problems and Theorems}

\newsiamthm{convariable}{Condition}
\newcommand{\bconvar}[1]{\renewcommand{\theconvariable}{#1}\setcounter{theorem}{\value{theorem}-1}\begin{convariable}}
\newcommand{\econvar}{\end{convariable}}
\newsiamthm{probvariable}{Problem}
\newcommand{\bprobvar}[1]{\renewcommand{\theprobvariable}{#1}\setcounter{theorem}{\value{theorem}-1}\begin{probvariable}}
\newcommand{\eprobvar}{\end{probvariable}}

\newcommand{\bconS}{\bconvar{S}}
\newcommand{\bconB}{\bconvar{B}}
\newcommand{\bconL}{\bconvar{L2}}
\newcommand{\bconAo}{\bconvar{A1}}
\newcommand{\bconLo}{\bconvar{L1}}
\newcommand{\bconA}{\bconvar{A2}}
\newcommand{\bconCo}{\bconvar{C1}}
\newcommand{\bconC}{\bconvar{C2}}
\newcommand{\bconAL}{\bconvar{AL}}
\newcommand{\bconK}{\bconvar{U}}

\newcommand{\bprobAE}{\bprobvar{AS}}
\newcommand{\bprobM}{\bprobvar{MAS}}
\newcommand{\bprobLT}{\bprobvar{SOAS}}
\newcommand{\bprobSVAR}{\bprobvar{SV}}

\newcommand{\bmps}{\begin{myprobsec}}
\newcommand{\emps}{\end{myprobsec}}
\newcommand{\bmp}{\begin{myprob}}
\newcommand{\emp}{\end{myprob}}


\begin{abstract}
We prove well-posedness results and a priori bounds on the solution of the Helmholtz equation $\nabla\cdot(A\grad u) + k^2 n u = -f$, posed either in $\RRd$ or in the exterior of a star-shaped Lipschitz obstacle, for a class of random $A$ and $n,$ random data $f$, and for all $k>0$. The particular class of $A$ and $n$ and the conditions on the obstacle ensure that the problem is nontrapping almost surely. These are the first well-posedness results and a priori bounds for the stochastic Helmholtz equation for arbitrarily large $k$ and for $A$ and $n$ varying independently of $k$. These results are obtained by combining recent bounds on the Helmholtz equation for deterministic $A$ and $n$ and general arguments (i.e.~not specific to the Helmholtz equation) presented in this paper for proving a priori bounds and well-posedness of variational formulations of linear elliptic stochastic PDEs. We emphasise that these general results do not rely on either the Lax-Milgram theorem or Fredholm theory, since neither are applicable to the stochastic variational formulation of the Helmholtz equation.
\end{abstract}



\section{Introduction}\label{sec:intro}
The goals of this paper are to prove results on the well-posedness of variational formulations of the stochastic Helmholtz equation
\beq\label{eq:hh-intro}
\grad\cdot\mleft(A(\omega)\grad u(\omega)\mright) + k^2n(\omega)u(\omega) = -f(\omega),
\eeq
as well as a priori bounds on its solution that are explicit in the wavenumber $k$ and the material coefficients $A$ and $n.$

We consider \eqref{eq:hh-intro} with physical domain either $\RRd, \, d=2,3,$ or $\RR^d\setminus\Dmclos,$ where $\Dm$ (referred to as the \defn{obstacle}) is a bounded, Lipschitz, open set such that $\RRd\setminus\Dmclos$ is connected, and

\bit
\item $\omega$ is an element of the underlying probability space,
\item $A$ is a symmetric-positive-definite matrix-valued random field such that $\supp(I-A)$ is compact,
\item $n$ is a positive real-valued random field such that $\supp(1-n)$ is compact,
\item $f$ is a real-valued random field such that $\supp f$ is compact, and
  \item $k>0$ is the wavenumber,
  \eit
and we are particularly interested in the case where the wavenumber $k$ is large.

\paragraph{Motivation} The motivation for establishing well-posedness and proving a priori bounds on the solution of \eqref{eq:hh-intro} is the growing interest in Uncertainty Quantification (UQ) for the Helmholtz equation; see e.g.~\cite{XiSh:07,TsXiYi:11,BuGh:14,GaHa:15,FeLiLo:15,FeLiNi:18,LiWaZh:18,HiScScSc:15,BaCaHaZh:18}. (In this PDE context, by `UQ' we mean theory and algorithms for computing statistics of quantities of interest involving PDEs \emph{either} posed on a random domain \emph{or} having random coefficients.) There is a large literature on UQ for the stationary diffusion equation
\beq\label{eq:diffusion}
-\grad\cdot (\kappa(\omega) \grad u(\omega))=f(\omega),
\eeq
due in part to its large number of applications (e.g.~in modelling groundwater flow), and a priori bounds on the solution are vital for the rigorous analysis of UQ algorithms; see e.g.~\cite{BaTeZo:04,BaNoTe:07,Gi:10,MuSt:11,ChScTe:13}. In contrast, whilst \eqref{eq:hh-intro} has many applications (e.g.~in geophysics and electromagnetics), there is much less rigorous theory of UQ for the Helmholtz equation. The main reason for this is that the (deterministic) PDE theory of \eqref{eq:hh-intro} when $k$ is large is much more complicated that the analogous theory for \eqref{eq:diffusion}.
 
 \paragraph{Related previous work} To our knowledge, the only work that considers \eqref{eq:hh-intro} with large $k$ and attempts to establish either (i) well-posedness of variational formulations or (ii) a priori bounds is \cite{FeLiLo:15}, which considers both (i) and (ii) for \eqref{eq:hh-intro} posed in a bounded domain with an impedance boundary condition. We discuss the results of \cite{FeLiLo:15} further in \cref{sec:otherwork}, but we highlight here that (a) \cite{FeLiLo:15} considers $A=I$ and $n=1+\eta,$ with $\eta$ random and the magnitude of $\eta$ decreasing with $k,$  whereas we consider classes of $A$ and $n$ that allow $k$-independent random perturbations, and (b) in its well-posedness result, \cite{FeLiLo:15}  invokes Fredholm theory to conclude existence of a solution, but this relies on an incorrect assumption about compact inclusion of Bochner spaces---see \cref{sec:federico} below. In \cref{sec:otherwork} we also discuss the papers \cite{BuGh:14,JeSc:16,JeScZe:17,HiScScSc:15} on the theory of UQ for either \eqref{eq:hh-intro} or the related time-harmonic Maxwell's equations; in these papers either the $k$-explicit well-posedness is not a primary concern or $k$ is assumed to be small. Our hope is that the results in the present paper can be used in the rigorous theory of UQ for Helmholtz problems with large $k.$
 
\paragraph{The contributions of this paper} The main results in this paper, \cref{thm:hh-gen,thm:hh-hetero} below, concern well-posedness and a priori bounds for the solutions of various formulations of the stochastic Helmholtz equation; these formulations include those used in sampling-based UQ algorithms (\cref{prob:msedp,prob:somsedp} below) and in the stochastic Galerkin method (\cref{prob:svsedp} below). These are the first such results for arbitrarily large $k$ and for $A$ and $n$ varying independently of $k$. These results are proved by combining:
\ben
\item bounds for the Helmholtz equation in \cite{GrPeSp:18} with $A$ and $n$ deterministic but spatially-varying, with
\item general arguments (i.e.~not specific to Helmholtz) presented here for proving a priori bounds and well-posedness of variational formulations of linear elliptic SPDEs.
\een
Regarding 1: the $k$-dependence of the bounds on $u$ in terms of $f$ depends crucially on whether or not $A$, $n$, and $\Dm$ are such that there exist trapped rays. In the trapping case, the solution operator can grow exponentially in $k$ (see \cite{Ra:71,Bu:98,PoVo:99,CaPo:02,Be:03a} and \cite[Section 2.5]{BeChGrLaLi:11}, and the reviews in \cite[Section 6]{MoSp:17}, \cite[Section 1.1]{ChSpGiSm:17}, and \cite[Section 1]{GrPeSp:18}); in contrast, in the nontrapping case, the solution operator is bounded uniformly in $k$ (see \cite{Va:75,MeSj:78,Bu:02}). The bounds in \cite{GrPeSp:18} are under conditions on $A,n,$ and $\Dm$ that ensure nontrapping of rays; the significance of these bounds is that they are the first (deterministic) bounds for the Helmholtz scattering problem in which both $A$ and $n$ vary and the bounds are explicit in $A$ and $n$ (as well as in $k$). This feature of being explicit in $A$ and $n$ is crucial in allowing us to prove the results in this paper when $A$ and $n$ are random fields.

Regarding 2: the main reason these general arguments are needed is the fact that the variational formulations of both the deterministic and the stochastic Helmholtz equation are not coercive, and so one cannot use the Lax--Milgram theorem to conclude well-posedness and an a priori bound.  In the deterministic case, the remedy for the lack of coercivity of the Helmholtz equation is to use Fredholm theory, but this is \emph{not} applicable to the stochastic variational formulation of the Helmholtz equation because the necessary compactness results do not hold in Bochner spaces (see \cref{sec:federico} below). Our solution to this lack of coercivity and failure of Fredholm theory is to use well-posedness results and bounds from the deterministic case to prove results for the stochastic case. We work `pathwise' by integrating the deterministic results over probability space, identifying conditions under which the necessary quantities are indeed integrable. Our approach is given in a general framework that, given (i) deterministic well-posedness results and a priori bounds that are explicit in all the coefficients, and (ii) measurability and integrability conditions on the stochastic quantities, returns corresponding well-posedness results, a priori bounds, and equivalence results for different formulations of the stochastic problem. One reason we state our well-posedness results in general (i.e.~not only in the specific case of the Helmholtz equation) is that we expect that they can be used in the future to prove well-posedness results for the time-harmonic Maxwell's equations in random media. A nontechnical summary of the ideas behind our general well-posedness results is given in \cref{rem:nontechnical} below. Some of these results are similar in spirit to the results about the PDE \eqref{eq:diffusion} in \cite{Gi:10,MuSt:11} (which deal with the failure of Lax--Milgram for the stochastic variational problem for \eqref{eq:diffusion} in the case when the coefficient $\kappa$ is not uniformly bounded above and below), and our general arguments use some of the ideas and technical tools from these two papers. 

\subsection{Statement of main results}\label{sec:hh-results}

\begin{figure}
\begin{centering}
\scalebox{0.65}{
\begin{tikzpicture}[scale=0.5,even odd rule]

    \pgfdeclarepatternformonly{owennortheast}
    {\pgfpoint{-0.1cm}{-0.1cm}}{\pgfpoint{1.1cm}{1.1cm}}
    {\pgfpoint{1cm}{1cm}}
    {
    \pgfpathmoveto{\pgfpointorigin}
    \pgfpathlineto{\pgfpoint{1cm}{1cm}}
    \pgfsetlinewidth{0.5\pgflinewidth}
    \pgfusepath{stroke}
    }
    
     \pgfdeclarepatternformonly{owennorthwest}
    {\pgfpoint{-0.1cm}{-0.1cm}}{\pgfpoint{1.1cm}{1.1cm}}
    {\pgfpoint{1cm}{1cm}}
    {
    \pgfpathmoveto{\pgfpoint{1cm}{0cm}}
    \pgfpathlineto{\pgfpoint{0cm}{1cm}}
   \pgfsetlinewidth{0.5\pgflinewidth}
    \pgfusepath{stroke}
    }
    
      \pgfdeclarepatternformonly{owennorthwestsmall}
    {\pgfpoint{-1cm}{-1cm}}{\pgfpoint{1cm}{1cm}}
    {\pgfpoint{0.5cm}{0.5cm}}
    {
    \pgfpathmoveto{\pgfpoint{1cm}{0cm}}
    \pgfpathlineto{\pgfpoint{0cm}{1cm}}
   \pgfsetlinewidth{0.5\pgflinewidth}
    \pgfusepath{stroke}
    }
    
       \pgfdeclarepatternformonly{owennortheastsmall}
    {\pgfpoint{-1cm}{-1cm}}{\pgfpoint{1cm}{1cm}}
    {\pgfpoint{0.5cm}{0.5cm}}
    {
    \pgfpathmoveto{\pgfpointorigin}
    \pgfpathlineto{\pgfpoint{1cm}{1cm}}
    \pgfsetlinewidth{0.5\pgflinewidth}
    \pgfusepath{stroke}
    }

    \filldraw[pattern=owennortheast] 
(4.5,0) .. controls (5,-1) and (4,-2) ..
(3,-2.5) .. controls (2,-3) and (2,-2.5) ..
(1.5,-3.5) .. controls (1,-4.5) and (-1,-3.5) ..
(-2,-3) .. controls (-3,-2.5) and (-2.5,-2.5) ..
(-4,-3) .. controls (-5.5,-3.5) and (-6,-0.5)..
(-5,0) .. controls (-4,0.5) and (-4,0.5) ..
(-4.5,1.5) .. controls (-5,2.5) and (-5,2.5) ..
(-3,3.5) .. controls (-1,2.5) and (1,3.5) ..
(1.5,3.5) .. controls (2,3.5) and (4,1) ..
    cycle;
    
    \draw (0,-1) circle [radius=10];
    
    \filldraw[pattern=owennorthwestsmall,yshift=-0.5cm]
    (0,8)  .. controls (1.5,9) and (2,9) ..
    (4,8) .. controls (5,7) and (4.5,7) ..
    (4,6.5) .. controls (3.5,6) and (4,6) .. 
    (4,5) .. controls (4,4) and (2,4) ..
    (0,5) .. controls (-2,6) and (-1.5,7) ..
    cycle;

    \filldraw[pattern=owennorthwestsmall]
    (1.5,-3.5) .. controls (1,-4.5) and (-1,-3.5) ..
    (-2,-3) .. controls (-3,-2.5) and (-2.5,-2.5) ..
    (-4,-3) .. controls (-5.5,-3.5) and (-8.8,-3)..
    (-9,-4) .. controls (-9.2,-5) and (-7,-5.5) ..
    (-6,-6) .. controls (-5,-6.5) and (-3,-6) ..
    (-2,-7) .. controls (-1,-8) and (1.5,-6) ..
    (2,-5.5) .. controls (2.5,-5) and (2,-4) ..
    cycle;
    
    \filldraw[pattern=owennortheastsmall]
    (3,-2.5) .. controls (2,-3) and (2,-2.5) ..
    (1.5,-3.5) .. controls (1,-4.5) and (-1,-3.5) ..
    (-2,-3) .. controls (-3,-2.5) and (1,-8) ..
    (2,-8) .. controls (3,-8) and (7,-8) ..
    (8,-6) .. controls (8.5,-5) and (8.5,-4) ..
    (8.5,-3.5) .. controls (7.5,-2.5) and (3.5,-2.5) ..
    cycle;


\draw (-0.1,0.1) node[fill=white] {\large$\Dm$};
\draw (7.8,6.8) node[fill=white] {\large$\GR$};
\draw (-5.3,-4.2) node[fill=white] {\large$\supp\mleft(I-A\mright)$};
\draw (5.3,-4.4) node[fill=white] {\large$\supp\mleft(1-n\mright)$};
\draw (1.5,5.8) node[fill=white] {\large$\supp f$};
\draw (6,2) node[fill=white] {\large$\DR$};
\end{tikzpicture}}
\caption{Examples of the domains $\Dm$  and $\DR$, the set $\GR,$ and essential supports of $I-A,$ $1-n$ and $f$ in the definition of the Helmholtz stochastic EDP.}  \label{fig:domain}
  \end{centering}
  \end{figure}
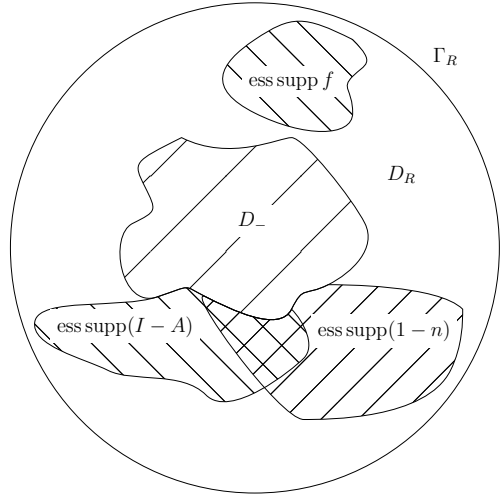

\paragraph{Notation and basic definitions}Let either (i) $\Dm \subset \RRd,$ $d=2,3,$ be a bounded Lipschitz open set such that $\bzero \in \Dm$ and the open complement $\Dp\de \RR^d\setminus \overline{\Dm}$ is connected, or (ii) $\Dm = \emptyset.$ Let $\GD = \partial \Dm.$ 
Fix $R>0$ and let $\BR$ be the ball of radius $R$ centred at the origin. Define $\GR := \partial \BR$ and $\DR \de \Dp \cap \BR$ (see \cref{fig:domain}). Let $\gamma$ denote the trace operator from $\DR$ to $\partial \DR = \GD \cup \GR$ and define $\HozDDR \de \set{v \in \HoDR \st \gamma v = 0 \ton \GD}.$ 
 
Let $\TrR: H^{1/2}(\Gamma_R) \rightarrow H^{-1/2}(\Gamma_R)$ be the Dirichlet-to-Neumann map for the deterministic equation $\Delta u+k^2 u=0$ posed in the exterior of $\BR$ with the Sommerfeld radiation condition 
\beq\label{eq:src}
\frac{\partial u}{\partial r}(\bx) - \ri ku(\bx) = o\mleft(\frac1{r^{(d-1)/2}}\mright) \text{ as } r\de\abs{\bx}\rightarrow \infty, \text{ uniformly in } \frac{\bx}{\abs{\bx}};
\eeq
see \cite[Section 3]{KeGi:89}, \cite[Section 2.6.3]{Ne:01}, and \cite[Equations 3.5 and 3.6]{ChMo:08} for an explicit expression for $\TrR$  in terms of Hankel functions and Fourier series ($d=2$)/spherical harmonics ($d=3$). Let $\IPGR{\cdot}{\cdot}$ be the duality pairing on $\GR$ between $\HmhGR$ and $\HhGR$ and write $\dd\Leb$ for Lebesgue measure.

We let $\spd$ be the set of all symmetric-positive-definite matrices in $\RRdtd$. Let $\LiDpRRdtd$ be the set of all matrix-valued functions $A:\Dp\rightarrow\spd$ such that $A_{i,j} \in \LiDpRR$ for all $i,j = 1,\ldots,d.$ Where the range of functions is $\CC$ we suppress the second argument in a function space, e.g.~we write $\LiDp$ for $\LiDpCC.$ For $\Az \in \spd,$ we write $\NopCCd{\Az}$ for the operator norm induced by the Euclidean vector norm on $\CCd$ (i.e., $\NopCCd{\cdot}$ is the spectral norm), and for $\Az:\DR \rightarrow \spd,$ we write $\NLiDRspd{\Az}$ for the norm $\NLiDRRR{\NopCCd{\Az(\bx)}}.$ Observe that $\NLiDRspd{\cdot}$ is a slight abuse of notation, as $\LiDRspd$ is not a vector space. Nevertheless $\NLiDRspd{\cdot}$ still fulfils the properties of a norm on $\LiDRspd.$ We define $\NWoiDRspd{\cdot}$ as the componentwise maximum of the $W^{1,\infty}$ norms of the components (c.f.~the norm equivalence in \cref{eq:normsdef} below).

We write $\Do \compcont \Dt$ if $\Do$ is a compact subset of the open set $\Dt.$ Let $\OFP$ be a complete probability space. Throughout this paper, unless stated otherwise we equip a topological space with its Borel $\sigma$-algebra. See \cref{app:mtBs} for a summary of the measure-theoretic concepts used in this paper. Let
\bit
\item $f:\Omega\rightarrow\LtDp$ be such that $\supp f \compcont \BR$ almost surely
\item $n:\Omega\rightarrow \LiDpRR$ be such that $\supp(1-n) \compcont \BR$ almost surely and there exist $\nmin, \nmax:\Omega\rightarrow \RR$ such that
$0 < \nmin(\omega) \leq n(\omega)(\bx) \leq \nmax(\omega)$
for almost every $\bx \in \Dp$ almost surely, and
\item $A:\Omega\rightarrow\LiDpRRdtd$ be such that $\supp(I-A) \compcont \BR,\, A_{ij} = A_{ji}$ almost surely, and there exist $\Amin,\Amax:\Omega \rightarrow \RR$ such that $0 < \Amin(\omega) < \Amax(\omega)$ almost surely and
$\Amin(\omega)\abs{\bxi}^2 \leq \big(A(\omega)(\bx)\bxi\big)\cdot\bxi \leq \Amax(\omega)\abs{\bxi}^2$
for almost every $\bx \in \Dp$ and for all $\bxi \in \CCd$ almost surely.
\eit
If $v:\Omega \rightarrow Z$ for some function space $Z$ of functions on $\RRd,$ we abuse notation slightly and write $v(\omega,\bx)$ instead of $v(\omega)(\bx).$

\paragraph{Variational Formulations} We consider three different formulations of the  \emph{Helmholtz stochastic exterior Dirichlet problem} (stochastic EDP); \cref{prob:msedp,prob:somsedp,prob:svsedp} below.

Define the sesquilinear form $a(\omega)$ on $\HozDDR \times \HozDDR$ by
\beq\label{eq:SEDPa}
\mleft[a(\omega)\mright]\mleft(\vo,\vt\mright)\de\int_{D_R}\Big( \mleft( A(\omega) \grad \vo\mright)\cdot \grad \vtb 
 - k^2 n(\omega)\, \vo\,\vtb \Big)\dd\Leb- \big\langle T_R \gamma \vo,\gamma \vt\big\rangle_{\Gamma_R},
 \eeq
 and the antilinear functional $L(\omega)$ on $\HozDDR$ by
\beq\label{eq:SEDPL}
\mleft[L(\omega)\mright](\vt)\de \int_{D_R} f(\omega)\, \vtb\,\dd\Leb.
\eeq
Define the sesquilinear form $\as$ on $L^2\big(\Omega;H_{0,D}^1(D_R)\big)\times L^2\big(\Omega;H_{0,D}^1(D_R)\big)$ and the antilinear functional $\Ls$ on $L^2\big(\Omega;H_{0,D}^1(D_R)\big)$ by 
\beq\label{eq:SEDPas}
\as\mleft(\vo,\vt\mright)\de \int_\Omega \mleft[a(\omega)\mright]\mleft(\vo(\omega),\vt(\omega)\mright)\dd\PP(\omega)
\quad\text{ and } \quad
%
\Ls(\vt)\de \int_\Omega \mleft[L(\omega)\mright]\mleft(\vt(\omega)\mright)\dd\PP(\omega).
\eeq
We consider the following three problems:

\bprobvar{1}[Measurable EDP almost surely]\label{prob:msedp}
Find a measurable $u:\Omega\rightarrow\HozDDR$ such that
\vspace{-2ex}
\beqs
\mleft[a(\omega)\mright]\mleft(u(\omega),v\mright) = \mleft[L(\omega)\mright](v) \tforall v \in \HozDDR \text{ almost surely.}
\eeqs
\eprobvar

\bprobvar{2}[Second-order EDP almost surely]\label{prob:somsedp}
Find $u\in L^2\big(\Omega;H_{0,D}^1(D_R)\big)$ such that
\beqs
\mleft[a(\omega)\mright]\mleft(u(\omega),v\mright) = \mleft[L(\omega)\mright](v) \tforall v \in \HozDDR \text{ almost surely.}
\eeqs
\eprobvar

\bprobvar{3}[Stochastic variational EDP]\label{prob:svsedp}
Find $u\in L^2\big(\Omega;H_{0,D}^1(D_R)\big)$ such that
\beqs
\as(u,v) = \Ls(v) \tforall v \in \LtOHozDDR.
\eeqs
\eprobvar

\Cref{prob:somsedp} is the foundation of sampling-based UQ methods, such as Monte-Carlo and Stochastic-Collocation methods; its analogue for the stationary diffusion equation is well-studied in, e.g., \cite{XiHe:05,BaNoTe:07,NoTeWe:08a,Ch:12,ChScTe:13,TeJaWeGu:15,KuNu:16,HeLaSc:18}. Similarly \cref{prob:svsedp} is the foundation of the Stochastic Galerkin method (a finite element method in $\Omega \times D,$ where $D$ is the spatial domain), and is studied for the Helmholtz Interior Impedance Problem in \cite{FeLiLo:15}, and its analogue for the stationary diffusion equation is considered in, e.g., \cite{BaTeZo:04,KhSc:11,BaScZo:11,GuWeZh:14}.

\bre[Why consider \cref{prob:msedp}?]\label{rem:whyone}

The difference between \cref{prob:msedp,prob:somsedp} is that \cref{prob:msedp} requires no integrability of $u$ over $\Omega$, whereas \cref{prob:somsedp} requires $u \in \LtOHozDDRnormal$. Since all the theory for sampling-based UQ methods assume some integrability of the solution, the natural question is: why consider \cref{prob:msedp} at all?

The main reason we consider \cref{prob:msedp} is that, given the existing PDE theory for the Helmholtz equation, we can prove existence of a solution to \cref{prob:msedp} under general conditions on $A$ and $n$, but there is no current prospect of proving existence of a solution to \cref{prob:somsedp} under general conditions on $A$ and $n$. The explanation for this consists of the following three points:
\ben
\item The only two known ways to obtain a solution to \cref{prob:somsedp} are: (i) obtain a deterministic a priori bound, explicit in all parameters, and integrate (followed, e.g., in \cite{ChScTe:13} for \eqref{eq:diffusion} with lognormal coefficients) and (ii) obtain a solution to \cref{prob:svsedp} and show this is a solution to \cref{prob:somsedp}. In the Helmholtz case, doing (ii) is difficult as neither the Lax--Milgram theorem nor Fredholm theory is applicable (as explained in the introduction), and so we follow the approach in (i).
\item The only known bounds on the solution of the Helmholtz equation explicit in all parameters are those recently obtained for nontrapping scenarios in \cite{GrPeSp:18,GaSpWu:18}.
\item Obtaining a bound explicit in all parameters for a general class of $A$ and $n$, e.g., $A \in \WoiDRspd$ and $n \in \LiDRRR$ is well beyond current techniques. Indeed, a general class of $A$ and $n$ will include both trapping and nontrapping scenarios, and such a bound would need to capture the exponential blow-up in $k$ for trapping $A$ and $n$, the uniform boundedness in $k$ for nontrapping $A$ and $n$, and be explicit in $A$ and $n$.
\een
Given this fact that there is no current prospect of proving existence of a solution to \cref{prob:somsedp} under general conditions on $A$ and $n$ we keep \cref{prob:msedp} so that we prove an (albeit weaker) existence result for the Helmholtz equation with general coefficients.
\ere

\bre[Measurability of $u$ in \cref{prob:msedp}]
It is natural to construct  the solution of \cref{prob:msedp} pathwise; that is, one defines $u(\omega)$ to be the solution of the deterministic problem with coefficients $A(\omega)$ and $n(\omega).$ However, it is then not obvious that $u$ is measurable.
In the proof of \cref{thm:hh-gen} below, we show that the measurability of $u$ follows from
\ben
\item a natural condition on the measurability of the coefficients and data (\cref{con:cborel} below), and 
\item the continuity of the map taking the coefficients of the deterministic PDE to the solution of the deterministic PDE (see \cref{lem:solcont} below).
  \een
\ere

In \cref{thm:hh-gen,thm:hh-hetero} we prove results on the well-posedness of \cref{prob:msedp,prob:somsedp,prob:svsedp} under conditions on $A,$ $n,$ $f,$ and $\Dm.$ Although $A,n,$ and $f$ are defined on $\Dp,$ since $\supp(I-A),$ $\supp(1-n),$ and $\supp f$ are compactly contained in  $\DR$ we can consider $A,n,$ and $f$ as functions on $\DR.$

\bcon[Regularity and stochastic regularity of $f,$ $A,$ and $n$]\label{con:hh-fAn}
The random fields $f, A,$ and $n$ satisfy $f \in \LtOLtDR,$  $A:\Omega \rightarrow \WoiDRspd$ with $A \in \LiOLiDRRRdtd,$ and $n \in \LiOLiDRRR.$
\econ

\bth[Equivalence of variational problems]\label{thm:hh-gen}
Under \cref{con:hh-fAn}:
\bit
\item The maps $\as$ and $\Ls$ (defined by \eqref{eq:SEDPas}) are well-defined.
\item $u \in
L^2\big(\Omega;H_{0,D}^1(D_R)\big)$
solves \cref{prob:somsedp} if and only if $u$ solves \cref{prob:svsedp}.
\item If $u \in
L^2\big(\Omega;H_{0,D}^1(D_R)\big)$
solves \cref{prob:somsedp}, then any member of the equivalence class of $u$ solves \cref{prob:msedp}.
\item The solution of \cref{prob:msedp} exists and is unique up to modification on a set of measure zero in $\Omega.$
\item The solution of \cref{prob:somsedp,prob:svsedp} is unique in $L^2\big(\Omega;H_{0,D}^1(D_R)\big)$.
\eit
\enth

Observe that the only relationship between formulations not proved in \cref{thm:hh-gen} is: if $u:\Omega\rightarrow \HozDDR$  solves \cref{prob:msedp} then $u \in 
L^2\big(\Omega;H_{0,D}^1(D_R)\big)$
and $u$ solves \cref{prob:somsedp}. \cref{thm:hh-hetero} below includes this relationship, but we need additional assumptions on $A,n,$ and $\Dm.$

\bde[A particular class of (deterministic) nontrapping coefficients]\label{def:hh-nontrapping}
Let $\muo,\mut >0,$ $\Az \in \WoiDRspd$ with $\supp(I-\Az) \compcont \BR$, and $\nz \in \WoiDRRR$ with $\supp(1-\nz) \compcont \BR.$ We write $\Az \in \NTA{\muo}$ and $\nz \in \NTn{\mut}$ if
\beq\label{eq:hh-Acond}
\Az(\bx) - \mleft(\bx\cdot\grad\mright)\Az(\bx) \geq \muo \quad \tand \quad \nz(\bx) + \bx\cdot\grad \nz(\bx) \geq \mut
\eeq
for almost every $\bx \in \DR,$ where the first inequality holds in the sense of quadratic forms.
\ede

The significance of the class of coefficients in \cref{def:hh-nontrapping} is that \cite[Theorem 2.5]{GrPeSp:18} proves bounds on the solution of \eqref{eq:hh-intro} for such $A$ and $n,$ where the constant in the bound only depends on $\muo, \mut, k, R, $ and $d.$

\bcon[{Nontrapping conditions on (random) $A$ and $n$}]\label{con:hh-hetero}
The random fields $A$ and $n$ satisfy $A:\Omega\rightarrow\WoiDRspd$ and $n:\Omega\rightarrow \WoiDRRR.$ Furthermore, there exist $\muo, \mut:\Omega\rightarrow \RR,$ independent of $f,$ with $\muo(\omega),\mut(\omega) > 0$ almost surely and $1/\muo,1/\mut \in \LtORR$  such that $A(\omega) \in \NTA{\muo(\omega)}$ almost surely and $n(\omega) \in \NTn{\mut(\omega)}$ almost surely.
\econ

\bde[Star-shaped]
The set $D \subseteq \RRd$ is \defn{star-shaped with respect to} the point $\bxz$ if for any $\bx \in D$ the line segment $\mleft[\bxz,\bx\mright] \subseteq D.$
\ede

\bth[{Equivalence of variational problems in a nontrapping case}]\label{thm:hh-hetero}
Let $\Dm$ be star-shaped with respect to the origin. Under \cref{con:hh-fAn,con:hh-hetero}:
\bit
\item The maps $\as$ and $\Ls$ (defined by \eqref{eq:SEDPas}) are well-defined.
\item \cref{prob:msedp,prob:somsedp,prob:svsedp} are all equivalent.
\item The solution $u \in
L^2\big(\Omega;H_{0,D}^1(D_R)\big)$
of these problems exists, is unique, and, given $\kz > 0,$ satisfies the bound
\beq\label{eq:Sbound1}
\NLtOLtDR{\grad u}^2 + k^2\NLtOLtDR{u}^2\leq \NLoO{\Co} \NLtOLtDR{f}^2
\eeq
for all $k\geq\kz$, where $\Co:\Omega\rightarrow\RR$ is given by
\beq\label{eq:C1}
\Co = \max\mleft\{\frac1{\mu_1},\frac1{\mu_2}\mright\}\mleft(\frac{R^2}{\mu_1} + \frac{2}{\mu_2}\mleft(R+ \frac{d-1}{2\kz}\mright)^2\mright).
\eeq
\eit
\enth

As highlighted above, \cref{thm:hh-hetero} is obtained from combining
deterministic a priori bounds from \cite{GrPeSp:18} with the general
arguments in \cref{sec:general} about well-posedness of variational
formulations of stochastic PDEs. \cref{thm:hh-hetero} uses the most basic a priori bound proved
in \cite{GrPeSp:18} (from \cite[Theorem 2.5]{GrPeSp:18}), but \cite{GrPeSp:18}
contains several extensions of this bound. \cref{rem:planewave,rem:ones,rem:tedp,rem:jumps,rem:kdep}
outline the implications that these (deterministic) extensions have for
the stochastic Helmholtz equation.

\bre[Dirichlet boundary conditions on $\GD$ and plane-wave incidence]\label{rem:planewave}
The formulations of the stochastic EDP above assume that $u=0$ on the boundary $\GD.$ An important scattering problem for which $u \neq 0$ on $\GD$ is when $u$ is the field scattered by an incident plane wave; in this case $\gamma u = -\gamma \uI,$ where $\uI$ is the incident plane wave \cite[p. 107]{ChGrLaSp:12}.

The results in this paper can be easily extended to the case when $u\neq0$ on $\pDm$ using \cite[Theorem 2.19(ii)]{GrPeSp:18}  which proves a priori (deterministic) bounds in this case. One subtlety, however, is that $f$ is then not necessarily independent of $\muo$ and $\mut.$ Indeed in this case
$f = -\grad\cdot\mleft(A\grad \uI\mright) - k^2n\uI$.
If $\muo$ depends on $A$ and $\mut$ depends on $n$ then
$f$ may be not be independent of $\muo$ and $\mut.$ One can produce an analogue of \cref{thm:hh-hetero} in the case where $f,\muo,$ and $\mut$ are dependent, but one requires $1/\muo, 1/\mut \in \LfO$ and $f \in \LfOLtD;$ see \cref{rem:notindep} below.
\ere

\bre[The case when either $n=1$ or $A=I$]\label{rem:ones}
When either $n=1$ or $A=I,$ \cite[Theorem 2.19]{GrPeSp:18} gives deterministic bounds under weaker conditions on $A$ and $n$ respectively; the corresponding results for the stochastic case are that:

\bit
\item 

When $n=1$  almost surely, the condition $A(\omega) \in \NTA{\muo(\omega)}$ in \cref{con:hh-hetero} can be improved to
$2A(\omega) - \mleft(\bx \cdot \nabla \mright)A(\omega) \geq \muo(\omega)$
for almost every $\bx\in \Dp,$ almost surely.

\item When $A=I$ almost surely, the condition $n(\omega) \in \NTn{\mut(\omega)}$ in \cref{con:hh-hetero} can be improved to:
\beq\label{eq:nimproved}
2n(\omega) + \bx \cdot \nabla n(\omega) \geq \mut(\omega) \,\text{ for almost every $\bx \in \Dp$, almost surely}.
\eeq
\eit
\ere

\bre[Geometric interpretation of the conditions on $A$ and $n$ in \cref{def:hh-nontrapping}]
Recall that the $k\tendi$ asymptotics of solutions of the Helmholtz
equation are governed by the behaviour of rays (see, e.g.,
\cite{BaBu:91}). Given (deterministic) $\Az$ and $\nz,$ the Helmholtz EDP is \defn{nontrapping} if all rays starting
in $\DR$ and evolving according to the
Hamiltonian flow defined by the symbol of $\grad\cdot\mleft(\Az\grad u\mright) + k^2\nz u = -\fz$  escape from $\DR$ after some uniform time (see, e.g., \cite[Definition 1.1]{Bu:02}); the EDP is \defn{trapping} otherwise.
The $k$-dependence of the solution operator depends strongly on
whether the problem is trapping, and the type of trapping present;
see, e.g., the overview discussions in \cite[Section 1]{GrPeSp:18},
\cite[Section 1.1]{ChSpGiSm:17}.

The conditions on $A$ and $n$ in \cref{con:hh-hetero} and the star-shapedness restriction
on $\Dm$ are sufficient for the Helmholtz stochastic EDP to be nontrapping almost surely. As noted in \cref{rem:ones}, when $A=I$ almost surely the condition on $n$  can be improved  from that in \eqref{eq:hh-Acond} to \eqref{eq:nimproved} using \cite[Theorem 2.19(ii)]{GrPeSp:18}.
The condition \eqref{eq:nimproved} is equivalent  to nontrapping when $n$ is radial, i.e. $n(\omega,\bx)=
n(\omega,\abs{\bx})$. Indeed, if $n$ is radial and $2n(\omega,\bx)+ \bx \cdot\nabla n(\omega,\bx)<0$ at a point $\bx \in \RRd$,
then the deterministic Helmholtz EDP given by $n(\omega,\bx)$ is trapping; see \cite{Ra:71} and \cite[Theorem 7.7]{GrPeSp:18}.
\ere

\bre[The Helmholtz stochastic truncated exterior Dirichlet problem]\label{rem:tedp}
When applying the Galerkin method to \cref{prob:msedp,prob:somsedp,prob:svsedp}, the Dirichlet-to-Neumann map $\TrR$ is expensive to compute. Therefore, it is
common to approximate the DtN map on $\GR$ by an `absorbing
boundary condition' (see, e.g., \cite[Section 3.3]{Ih:98} and the references
therein), the  simplest of which is the impedance boundary condition
$\partial u/\partial \nu - \ri k u=0$. We call the Helmholtz stochastic EDP posed in $\DR$ with
an impedance boundary condition on $\GR$ the stochastic \emph{truncated
exterior Dirichlet problem} (stochastic TEDP). In fact, since we no longer need to know
the DtN map explicitly on the truncation boundary, the truncation
boundary can be arbitrary (i.e. it does not have to be just a circle/sphere). Note that in the case when the obstacle is the empty set, the TEDP is just the Interior Impedance Problem.

The results in this paper also hold for the stochastic TEDP (with arbitrary Lipschitz truncation boundary) under an analogue of \cref{con:hh-hetero} based on the deterministic bounds in \cite[Theorem A.6(i)]{GrPeSp:18} instead of \cite[Theorem 2.5]{GrPeSp:18}.
\ere

\bre[Discontinuous $A$ and $n$]\label{rem:jumps}
The requirements on $A$ and $n$ in \cref{def:hh-nontrapping} require them to be continuous (since  $\WoiDR = \CzoDR$ as $\DR$ is Lipschitz; see, e.g., \cite[Section 4.2.3, Theorem 5]{EvGa:92}). In addition to proving deterministic a priori bounds for the class of $A$ and $n$ in \cref{def:hh-nontrapping}, the paper \cite{GrPeSp:18} proves deterministic bounds for discontinuous $A$ and $n$ satisfying \eqref{eq:hh-Acond} in a distributional sense; see \cite[Theorem 2.7]{GrPeSp:18}. In this case, when moving outward from the obstacle to infinity,  $A$ can jump downwards and $n$ can
jump upwards on interfaces that are star-shaped. (When the jumps are in the opposite direction, the problem is trapping; see \cite{PoVo:99} and \cite[Section 6]{MoSp:17}). The
well-posedness results and a priori bounds in this paper can therefore be adapted to prove results about the stochastic Helmholtz equation
for a class of random $A$ and $n$ that allows nontrapping jumps on randomly-placed star-shaped interfaces.
\ere

\bre[$k$-dependent $A$ and $n$]\label{rem:kdep}
In this paper we focus on random fields $A$ and $n$ varying independently of $k;$ this corresponds to a fixed physical
medium, characterised by $A$ and $n$, with waves of frequency $k$ passing through.
In \cref{sec:generating} below we construct $A$ and $n$ as ($k$-independent) $W^{1,\infty}$ perturbations of random fields $\Az$ and $\nz$ satisfying \cref{con:hh-hetero}.
We note, however, that results for $A$ and $n$ being
\emph{$k$-dependent} $L^\infty$ perturbations (i.e. rougher, but $k$-dependent perturbations) of $\Az$
and $\nz$ satisfying \cref{con:hh-hetero} can easily be obtained.

The basis for these bounds is observing that \emph{deterministic} a priori bounds hold when
(a) $A\in \NTA{\mu_1}$, $n = \nz + \eta,$ where $\nz \in \NTn{\mu_2}$
and $k\NLiDRRR{\eta}$ is sufficiently small, and
(b) $A=\Az+B$, $n=\nz+\eta$, where $\Az\in \NTA{\mu_1}$, $\nz \in
\NTn{\mu_2}$, $k\NLiDRRR{\eta}$ and $k \NWoiDRspd{B}$ are
both sufficiently small, and $A, n,$ and $D_-$ are such that $u\in
H^2(D_R)$ (see, e.g., \cite[Theorem 4.18(i)]{Mc:00} or \cite[Theorems
2.3.3.2 and 2.4.2.5]{Gr:85} for these latter requirements). Given these deterministic bounds, the general arguments in this paper
can then be used to prove well-posedness of the analogous stochastic
problems.

To understand why bounds hold in the case (a), observe that one can
write the PDE as
\beq\label{eq:pert}
\nabla \cdot(A\grad u) +k^2 n_0 u = -f - k^2 \eta u;
\eeq
if $k\NLiDRRR{\eta}$ is sufficiently small then the contribution from
the $k^2 \eta u$ term on the right-hand side of \eqref{eq:pert} can be
absorbed into the $k^2\|u\|^2_{L^2(D_R)}$ term appearing on the left-hand
side of the bound (the deterministic analogue of \eqref{eq:Sbound1}).
In the case $\nz=1,$ this is essentially the argument used to prove the
a priori bound in \cite[Theorem 2.4]{FeLiLo:15} (see \cite[Remark
2.15]{GrPeSp:18}).
The reason bounds hold in the case (b) is similar, except now we need
the $H^2$ norm of $u$ on the left-hand side of the bound (as well as
the $H^1$ norm) to absorb the contribution from the $\nabla \cdot
(B\grad u)$ term on the right-hand side.
\ere

\subsection{Random fields satisfying \cref{con:hh-hetero}}\label{sec:generating}
The main focus of this paper is proving well-posedness of the variational formulations of the stochastic Helmholtz equation, and a priori bounds on the solution, for the most-general class of $A$ and $n$ allowed by the deterministic bounds in \cite{GrPeSp:18}. However, in this section, motivated by the Karhunen-Lo\`eve expansion (see e.g.~\cite[p.~201ff.]{LoPoSh:14}) and similar expansions of material coefficients for the stationary diffusion equation \cite[Section 2.1]{KuNu:16}, we consider $A$ and $n$ as series expansions around known non-random fields $\Az$ and $\nz$ satisfying \cref{con:hh-hetero} (i.e., \cref{con:hh-hetero} is satisfied for $\nz, \Az$ independent of $\omega \in \Omega$, and therefore $\muo,\mut$ independent of $\omega$).
Define
\beq\label{eq:nseries}
A(\omega,\bx) = \Az(\bx) + \sum_{j=1}^\infty \Yj(\omega) \Psij(\bx)\quad\text{and}\quad n(\omega,\bx) = \nz(\bx) + \sum_{j=1}^\infty\Zj(\omega) \psij(\bx) ,
\eeq
where:
\bit
\item $\supp\mleft(1-\Az\mright),\,\supp \mleft(I-\nz\mright) \compcont \BR,$
\item $\Az$ and $\nz$ satisfy \cref{con:hh-hetero} with $\muo$ and $\mut$ independent of $\omega \in \Omega$
\item $\Yj,\Zj \sim \Unif(-1/2,1/2)$ i.i.d.,
\item $\Psij \in \WoiDRspd$ with $\supp \Psij \compcont \BR$ for all $j =1,\ldots,\infty$,
\beq\label{eq:Apsimeas}
\sum_{j=1}^\infty \NWoiDRRRdtd{\Psij} < \infty, \tand
\eeq
\beq\label{eq:Apsipos}
\sum_{j=1}^\infty \esssup_{\bx \in \DR} \NopCCd{\Psij} < 2\Azmin,
\eeq
where $\Azmin > 0$ is such that $\Azmin\abs{\bxi}^2 \leq \big(A(\bx)\bxi\big)\cdot \bxi$ for almost every $\bx \in \Dp$ and for all $\bxi \in \CCd$, and where $\NopCCd{\cdot}$ is the operator norm induced by the Euclidean vector norm on $\CCd$ (i.e., $\NopCCd{\cdot}$ is the spectral norm).
\item $\psij \in \WoiDRRR$ with $\supp \psij \compcont \BR$ for all $j = 1,\ldots,\infty$,
\beq\label{eq:npsimeas}
\sum_{j=1}^\infty \NWoiDRRR{\psij} < \infty, \tand
\eeq
\beq\label{eq:npsipos}
\sum_{j=1}^\infty \NLiDRRR{\psij} < 2 \nzmin,
\eeq
where $\nzmin \de \essinf_{\bx \in \DR} \nz(\bx),$ and
\eit

The assumptions \eqref{eq:Apsipos} and \eqref{eq:npsipos} ensure that $A > 0$ (in the sense of quadratic forms) and $n > 0$  almost surely, and the assumptions \eqref{eq:Apsimeas} and \eqref{eq:npsimeas} are used to prove $A$ and $n$ are measurable.

Regarding the measurability of $A$ and $n$ defined by \eqref{eq:nseries}: the proof that $A$ and $n$ given by \eqref{eq:nseries} are measurable is given in \cref{lem:seriesmeas}, and relies on the proof that the sum of measurable functions
is measurable. This latter result is standard, but we have not been able to find this result for this particular setting of mappings
into a separable subspace of a general normed vector space, and so we briefly give it in \cref{lem:summeas}. 

The following lemmas give sufficient conditions for the series in \eqref{eq:nseries} to satisfy \cref{con:hh-hetero}.

\ble[Series expansion of $A$ satisfies \cref{con:hh-hetero}]\label{lem:Agen}
Let $\mu > 0$, $\delta \in \mleft(0,1\mright).$ If $\Az \in \NTA{\mu},$ and
\beq\label{eq:Aseriescond}
\sum_{j=1}^\infty \esssup_{\bx \in \DR} \NopCCd{\Psij(\bx) - \mleft(\bx\cdot\grad\mright)\Psij(\bx)} \leq 2\delta\mu,
\eeq
then $A \in \NTA{(1-\delta)\mu}$ almost surely.
\ele

\bpf[Proof of \cref{lem:Agen}]
Since $\Az \in \NTA{\mu},$ we have 
\beq\label{eq:Aseries1}
\Big(\mleft(A(\omega,\bx) - \mleft(\bx\cdot\grad\mright)A(\omega,\bx)\mright)\bxi\Big)\cdot\bxib \geq\mu \abs{\bxi}^2+\sum_{j=1}^\infty  \Big(\Yj(\omega)\mleft(\Psij(\bx)- \mleft(\bx\cdot\grad\mright)\Psij(\bx)\mright)\bxi\Big)\cdot\bxib
\eeq
for all $\bxi \in \CCd,$ for almost every $\bx \in \DR,$ almost surely. As $\Yj \sim \Unif(-1/2,1/2)$ for all $j$ and the bound \eqref{eq:Aseriescond} holds, the right-hand side of \eqref{eq:Aseries1} is bounded below by
\beqs
\mu \abs{\bxi}^2 -\half 2\delta \mu \abs{\bxi}^2 = (1-\delta)\mu\abs{\bxi}^2\text{ almost surely.}
\eeqs
Since $\bxi \in \CCd$ was arbitrary, it follows that $A(\omega) \in \NTA{(1-\delta)\mu)}$ almost surely, as required.
\epf

\ble[Series expansion of $n$ satisfies \cref{con:hh-hetero}]\label{lem:ngen}
Let $\mu > 0$ and $\delta \in \mleft(0,1\mright).$ If $\nz \in \NTn{\mu}$ and
\beq\label{eq:nseriescond}
\sum_{j=1}^\infty\NLiDRRR{\psij(\bx) + \bx\cdot\grad\psij(\bx)} \leq 2\delta\mu,
\eeq
then $n \in \NTn{(1-\delta)\mu}.$
\ele

The proof of \cref{lem:ngen} is omitted, since it is similar to the proof of \cref{lem:Agen}; in fact it is simpler, because it involves scalars rather than matrices.

\subsection{Discussion of the main results in the context of other work on UQ for time-harmonic wave equations}\label{sec:otherwork}

In this section we discuss existing results on well-posedness of \eqref{eq:hh-intro}, as well as analogous results for the elastic wave equation and the time-harmonic Maxwell's equations. The most closely-related work to the current paper is \cite{FeLiLo:15} (and its analogue for elastic waves \cite{FeLo:17}), in that a large component of \cite{FeLiLo:15} consists of attempting to prove well-posedness and a priori bounds for the stochastic variational formulation (i.e.~\cref{prob:svsedp}) of the Helmholtz Interior Impedance Problem; i.e., \eqref{eq:hh-intro} with $A=I$ and stochastic $n$ posed in a bounded domain with an impedance boundary condition $\partial u/\partial \nu - ik u = g$ (recall that this boundary condition is a simple approximation to the Dirichlet-to-Neumann map $\TrR$ defined above \eqref{eq:src}). Under the assumption of existence, \cite{FeLiLo:15} shows that for any $k>0$ the solution is unique and satisfies an a priori bound of the form \eqref{eq:Sbound1} (with different constant $\Co$), provided $n=1+\eta$ where the random field $\eta$ satisfies (almost surely) $\N{\eta}_{L^\infty} \leq C/k$ for some $C>0$ independent of $k$. \cite{FeLiLo:15} then invokes Fredholm theory to conclude existence, but this relies on an incorrect assumption about compact inclusion of Bochner spaces---see \cref{sec:federico} below. However, combining \cref{thm:hh-gen,rem:tedp,rem:kdep} with $A=I$ and $\nz=1+\eta$ (with $\eta$ as above) produces an analogous result to \cref{thm:hh-hetero}, and gives a correct proof of \cite[Theorem 2.5]{FeLiLo:15}. Therefore the analysis of the Monte Carlo interior penalty discontinuous Galerkin method in \cite{FeLiLo:15} can proceed under the assumptions of \cref{thm:hh-gen,rem:tedp,rem:kdep}.

The paper \cite{HiScScSc:15} considers the Helmholtz transmission problem with a stochastic interface, i.e.~\eqref{eq:hh-intro} posed in $\RRd$ with both $A$ and $n$ piecewise constant and jumping on a common, randomly-located interface. A component of this work is establishing well-posedness of \cref{prob:msedp} for this setup. To do this, the authors make the assumption that $k$ is small (to avoid problems with trapping mentioned above---see the comments after \cite[Theorem 4.3]{HiScScSc:15}); the sesquilinear form $a$ is then coercive and an a priori bound (in principle explicit in $A$ and $n$) follows \cite[Lemma 4.5]{HiScScSc:15}. By \cref{rem:jumps}, the results of this paper can be used to obtain the analogous well-posedness result for large $k$ in the case of nontrapping jumps.

The paper \cite{BuGh:14} studies the \defn{Bayesian inverse problem} associated to \eqref{eq:hh-intro} with $A=I$ and $n=1$  posed in the  exterior of a Dirichlet obstacle. That is, \cite{BuGh:14} analyses computing the posterior distribution of the shape of the obstacle given noisy observations of the acoustic field in the exterior of the obstacle. A component of the analysis in \cite{BuGh:14} is the well-posedness of the forward problem for an obstacle with a variable boundary \cite[Proposition 3.5]{BuGh:14}. Instead of mapping the problem to one with  a fixed domain and variable $A$ and $n,$ \cite{BuGh:14} works with the variability of the obstacle directly, using boundary-integral equations. The $k$-dependence of the solution operator is not considered, but would enter in \cite[Lemma 3.1]{BuGh:14}.

The papers \cite{JeScZe:17} and \cite{JeSc:16} consider the time-harmonic Maxwell's equations with (i) the material coefficients $\eps,\mu$ constant in the exterior of a perfectly-conducting random obstacle and (ii) $\eps,\mu$ piecewise constant and jumping on a common randomly located interface; in both cases these problems are mapped to problems where the domain/interface is fixed and $\eps$ and $\mu$ are random and heterogeneous. The papers \cite{JeScZe:17} and \cite{JeSc:16} essentially consider the analogue of \cref{prob:msedp} for the time-harmonic Maxwell's equations, obtaining well-posedness from the corresponding results for the related deterministic problems.

\subsection{Outline of the paper} In \cref{sec:otherwork} we discuss our results in the context of related literature. In \cref{sec:general} we  state general results on a priori bounds and well-posedness for stochastic variational formulations. In \cref{sec:genproof} we prove the results in \cref{sec:general}. In \cref{sec:hhproof} we prove \cref{thm:hh-gen,thm:hh-hetero}. In \cref{sec:federico} we discuss the failure of Fredholm theory for the stochastic variational formulation of Helmholtz problems. In \cref{app:mtBs} we recap results from measure theory and the theory of Bochner spaces.

\section{General results on proving a priori bounds and well-posedness of stochastic variational formulations}\label{sec:gen-framework}\label{sec:general}
In this section we state general results for proving a priori bounds and well-posedness results for variational formulations of linear elliptic SPDEs.

\subsection{Notation and definitions of the variational formulations}\label{sec:notdef}
Let $\OFP$ be a complete probability space. Let $X$ and $Y$ be separable Banach spaces over a field $\FF,$ (where $\FF = \RR$ or $\CC$).
Let $\homspace$ denote the space of bounded linear maps $X\rightarrow\Ys.$ Let $\cC$ be a topological space with topology $\TC.$ Given maps
\beqs
\coeff:\Omega\rightarrow\cC,\quad\toform:\cC \rightarrow \homspace,\quad\text{and } \rhs:\cC \rightarrow \Ys,
\eeqs let $\SAc:\LtOX \rightarrow \LtOYas$ and  $\SLc \in \LtOYas$ be defined by
\beq\label{eq:SA}
\big[\SAc(u)\big](v) \de \int_\Omega \big[\Acomega u(\omega)\big]\big(v(\omega)\big) \dd\PP(\omega) \quad\text{and} \quad 
\SLc(v) \de \int_\Omega \Lcomega\big(v(\omega)\big) \dd \PP(\omega)
\eeq
for $v \in \LtOY.$ Recall that a bounded linear map $X\rightarrow\Ys$ is equivalent to a sesquilinear (or bilinear) form on $X \times Y;$ see e.g. \cite[Lemma 2.1.38]{SaSc:11}. To keep notation compact, we write $\Acomega=(\toform \circ c)(\omega)$ and $\Lcomega=(\rhs\circ c)(\omega).$

\bre[Interpretation of the space $\cC$]
The space $\cC$ is the `space of inputs'. For the stochastic Helmholtz EDP in \cref{sec:hh-results} the space $\cC$ is defined in \cref{def:cCHh} below, but the upshot of this definition is that for any $\omega \in \Omega$ the triple $(A(\omega),n(\omega),f(\omega))$ is an element of $\cC.$
The maps $\coeff,$ $\toform,$ and $\rhs$ are given by $\coeff = \mleft(A,n,f\mright),$ $\toform = a,$ and $\rhs = L,$ where $a$ and $L$ are given by \eqref{eq:SEDPa} and \eqref{eq:SEDPL} respectively and the equality $\toform=a$ is meant in the sense of the one-to-one correspondence between $\homspace$ and sesquilinear forms on $X\times Y.$
\ere

The following three problems are the analogues in this general setting of \cref{prob:msedp,prob:somsedp,prob:svsedp} in \cref{sec:intro}.

\bprobM[Measurable variational formulation almost surely]\label{prob:meas}
Find a measurable function $u:\Omega \rightarrow X$ such that 
\beq\label{eq:aeeq}
\Acomega u(\omega) = \Lcomega \tin \Ys
\eeq
almost surely.
\eprobvar

\bprobLT[Second-order moment variational formulation almost surely]\label{prob:lt}
Find $u \in \LtOX$ such that \eqref{eq:aeeq} holds almost surely.
\eprobvar

\bprobSVAR[Stochastic variational formulation]\label{prob:svar}
Find $u \in \LtOX$ such that
\beq\label{eq:stoeq}
\SAc u = \SLc \tin \LtOYas.
\eeq
\eprobvar

\bre[Immediate relationships between formulations]\label{rem:imm}
Since $\LtOX \subseteq \BorelOX$ (the space of all measurable functions $\Omega\rightarrow X$)
it is immediate that if $u$ solves \cref{prob:lt} then every member of the equivalence class of $u$ solves \cref{prob:meas}.
\ere

\subsection{Conditions on $\toform,$ $\rhs,$ and $\coeff$}\label{sec:cons}
We now state the conditions under which we prove results about the equivalence of \cref{prob:meas,prob:lt,prob:svar}.

\bconAo[$\toform$ is continuous] \label{con:coeffstoform}
The function $\toform:\cC \rightarrow \homspace$ is continuous, where we place the norm topology on $X,$ the dual norm topology on $\Ys$, and the operator norm topology on $\homspace.$
\econvar

\bconA[Regularity of {$\toform\circ\coeff$}]\label{con:A}
The map $\toform\circ\coeff \in \LiOhomspace.$
\econvar

We note that \cref{con:A} is violated in the well-studied case of a log-normal coefficient $\kappa$ for the stationary diffusion equation \eqref{eq:diffusion}; in order to ensure the stochastic variational formulation is well-defined in this case, one must change the space of test functions as in \cite{Gi:10,MuSt:11}

\bconLo[$\rhs$ is continuous] \label{con:coeffstofunc}
The function $\rhs:\cC \rightarrow \Ys$  is continuous, where we place the dual norm topology on $\Ys.$
\econvar

\bconL[Regularity of {$\rhs\circ\coeff$}]\label{con:L}
The map $\rhs\circ \coeff \in \LtOYs.$
\econvar

\bconCo[$\coeff$ is measurable]\label{con:cborel}
 The function $c:\Omega \rightarrow \cC$ is measurable.
\econvar

To state the next condition, we need to recall the following definition.

\bde[{$\PP$}-essentially separably valued {\cite[p26]{Ry:02}}]\label{def:sepval}
Let $\mleft(S,\Top{S}\mright)$ be a topological space. A function $h:\Omega\rightarrow S$ is \defn{$\PP$-essentially separably valued} if there exists $E \in \cF$ such that $\PP(E) = 1$ and $h(E)$ is contained in a separable subset of $S$.
\ede

\bconC[{$\coeff$} is {$\PP$}-essentially separably valued]\label{con:C}
The map $\coeff:\Omega\rightarrow\cC$  is $\PP$-essentially separably valued.
\econvar

\bre[Why do we need \cref{con:C}?]
The theory of Bochner spaces requires strong measurability of functions (see \cref{def:strongmeas,def:bochnerspace} below). However, the proof techniques used in this paper rely heavily on the measurability of functions (see \cref{def:meas} below). In separable spaces these two notions are equivalent (see \cref{cor:pettis}). However, some of the spaces we encounter (such as $\LiDRRR$) are not separable. Therefore, in our arguments we use \cref{con:C} along with the Pettis Measurability Theorem (\cref{thm:pettis} below) to conclude that measurable functions are strongly measurable.
\ere

\bconB[A priori bound almost surely]\label{con:B}
There exist $\Cj,\fj:\Omega\rightarrow \RR, \,j=1,\ldots,m$ such that $\Cj\fj \in \LoO$ for all $j=1,\ldots,m$ and the bound
\beq \label{eq:sbe1}
\NX{u(\omega)}^2 \leq \sum_{j=1}^m\Cj(\omega)\fj(\omega)
\eeq
holds almost surely.
\econvar

\bre[Notation in the a priori bound]
We use the notation $f_j$ in the right-hand side of \eqref{eq:sbe1} to emphasise the fact that typically these terms relate to the right-hand sides of the PDE in question. For the stochastic Helmholtz EDP, $m=1,$ $f_1 = \NLtD{f}^2,$ and $C_1$ is given by \eqref{eq:C1}.
\ere

\bconK[Uniqueness almost surely]\label{con:K}
$\ker\mleft(\Acomega\mright) = \set{0}$ $\PP$-almost surely.
\econvar

The condition  $\ker\mleft(\Acomega\mright) = \set{0}$ $\PP$-almost surely can be stated as: given $\cG \in \LtOYas,$ for $\PP$-almost every $\omega \in \Omega$ the deterministic problem $\Acomega \uz = \cG$ has a unique solution,

\subsection{Results on the equivalence of \cref{prob:meas,,prob:lt,,prob:svar}}

\bth[Measurable solution implies second-order solution]\label{thm:3}
Under \cref{con:B}, if $u$ solves \cref{prob:meas} then $u$ solves \cref{prob:lt}  and satisfies the stochastic a priori bound
\beq\label{eq:sbresult}
\NLtOX{u}^2 \leq\sum_{j=1}^m \NLoO{C_jf_j}.
\eeq
\enth

Note that the right-hand side of the stochastic a priori bound \eqref{eq:sbresult} is the expectation of the right-hand side of the bound \eqref{eq:sbe1}.

\ble[Stochastic variational formulation well-defined]\label{lem:svarwelldefined}
Under \cref{con:coeffstoform,,con:A,,con:coeffstofunc,,con:L,,con:cborel,,con:C},  the maps $\SA$ and $\SL$ defined by \eqref{eq:SA} are well-defined in the sense that
\beq\label{eq:finite}
\mleft[\SA(\vo)\mright](\vt),\, \SL(\vt) < \infty \quad\text{for all } \vo \in \LtOX, \text{ for all }\vt \in \LtOY.
\eeq
\ele

\bth[Second-order solution implies stochastic variational solution]\label{thm:11}
Under \cref{con:coeffstofunc,,con:L,,con:cborel,,con:C}, if $u$ solves \cref{prob:lt} then $u$ solves \cref{prob:svar}.
\enth

\bth[Stochastic variational solution implies second-order solution]\label{thm:12}
If \cref{prob:svar} is well-defined and $u$ solves \cref{prob:svar}, then $u$ solves \cref{prob:lt}.
\enth

\begin{figure}[h]
  \centering
  \scalebox{0.90}{
\begin{tikzpicture}


\draw (0,1) node [rounded rectangle, fill=gray!45!white] (meas) {\Cref{prob:meas}};

\draw (0,-1) node [rounded rectangle, fill=gray!45!white] (lt) {\Cref{prob:lt}};

\draw (0,-3) node [rounded rectangle, fill=gray!45!white] (svar) {\Cref{prob:svar}};

\path node [left = \owenshift of meas.north] (meas top left) {};
\path node [right = \owenshift of meas.north] (meas top right) {};
\path node [left = \owenshift of meas.south] (meas bottom left) {};
\path node [right = \owenshift of meas.south] (meas bottom right) {};

\path node [left = \owenshift of lt.north] (lt top left) {};
\path node [right = \owenshift of lt.north] (lt top right) {};
\path node [left = \owenshift of lt.south] (lt bottom left) {};
\path node [right = \owenshift of lt.south] (lt bottom right) {};

\path node [left = \owenshift of svar.north] (svar top left) {};
\path node [right = \owenshift of svar.north] (svar top right) {};

\begin{scope}[->]

\draw
(meas bottom right)
--
node[right,align=center,text width=5cm] {Under \cref{con:B}, get stochastic a priori bound \eqref{eq:sbresult} (\Cref{thm:3})}
 (lt top right);

\draw
(lt top left)
--
node[left] {Immediate}
 (meas bottom left);

\draw
(lt bottom right)
--
node[right,align=center,text width=5cm] {Under \cref{con:coeffstofunc,,con:L,,con:cborel,,con:C},  (\Cref{thm:11})}
 (svar top right);

\draw
(svar top left)
--
node[left,align=center,text width=5cm] {If \cref{prob:svar} is well-defined (\cref{thm:12})}
(lt bottom left);
\end{scope}

\path node [align=center,below = \owentextshift of svar,text width=10cm] (svar wd) {Well-defined under  \cref{con:coeffstoform,,con:A,,con:coeffstofunc,,con:L,,con:cborel,,con:C} (\Cref{lem:svarwelldefined})};

\end{tikzpicture}
}
\caption{The relationship between the variational formulations. An arrow from Problem P to Problem Q with Conditions R indicates `under Conditions R, the solution of Problem P is a solution of Problem Q'}\label{fig:ladder}
\end{figure}
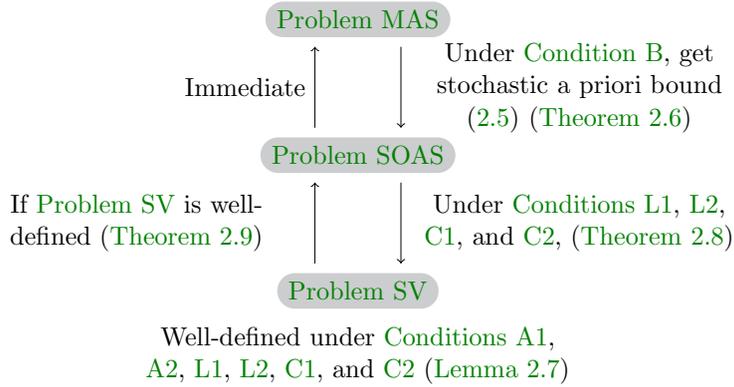

\Cref{thm:3,thm:11,thm:12,lem:svarwelldefined} are summarised in \cref{fig:ladder}.

\bre[\Cref{con:L} in \cref{thm:11}]
In \Cref{thm:11} we could replace \cref{con:L} with \cref{con:A}, and the result would still hold---see the proof for further details. However, \cref{con:L} is less restrictive than \cref{con:A}, as it only requires $L^2$ integrability of $\rhs\circ\coeff$ as opposed to essential boundedness of $\toform\circ\coeff.$ 
\ere

\ble[Showing uniqueness of the solution to \cref{prob:meas,prob:lt,prob:svar}]\label{lem:uniq}
If \cref{con:K} holds, then
\ben
\item\label[part]{part:uniq1} the solution to \cref{prob:meas} (if it exists) is unique up to modification on a set of $\PP$-measure 0 in $\Omega$, 
\item\label[part]{part:uniq2} the solution to \cref{prob:lt} (if it exists) is unique in $\LtOX,$ and 
\item\label[part]{part:uniq3} if \cref{prob:svar} is well-defined, the solution to \cref{prob:svar} (if it exists) is unique in $\LtOX.$
\een
\ele

\bre[Informal discussion on the ideas behind the equivalence results]\label{rem:nontechnical}
The diagram in \cref{fig:ladder} summarises the relationships between the variational formulations, and the conditions under which they hold. Moving `up' the left-hand side of the diagram, we prove a solution of \cref{prob:svar} is a solution of \cref{prob:lt} in \cref{thm:12}; the key idea in this theorem is to use a particular set of test functions and the general measure-theory result of \cref{lem:gotoae} below; this approach was used for the stationary diffusion equation \eqref{eq:diffusion} with log-normal coefficients in \cite{Gi:10}, and for a wider class of coefficients in \cite{MuSt:11}.

Moving `down' the right-hand side, we prove a solution of \cref{prob:meas} is a solution of \cref{prob:lt} in \cref{thm:3}; the key part of this proof is that the bound in \cref{con:B} gives information on the integrability of the solution $u.$ (In the case of \eqref{eq:diffusion} with uniformly coercive and bounded coefficient $\kappa,$ the analogous integrability result follows from the Lax--Milgram theorem; \cite[Proposition 2.4]{Ch:12} proves an equivalent result for \eqref{eq:diffusion} with lognormal coefficient $\kappa$ with an isotropic Lipschitz covariance function.) Proving a solution of \cref{prob:lt} is a solution of \cref{prob:svar} in \cref{thm:11} essentially amounts to posing conditions such that the quantities $\mleft[\Acomega\mleft(u(\omega)\mright)\mright] \mleft(v(\omega)\mright)$ and $\Lcomega \mleft(v(\omega)\mright)$ are Bochner integrable for any $v \in \LtOY,$ so that \eqref{eq:stoeq} makes sense. \Cref{lem:svarwelldefined} shows that the stronger property \eqref{eq:finite} holds, and requires stronger assumptions than \cref{thm:11}, since the proof of \cref{thm:11} uses the additional information that u solves \cref{prob:lt}.
\ere

\bre[Changing the condition $u \in \LtOX$]
Here we seek the solution $u \in \LtOX$ but we could instead require $u \in \LpOX,$ for some $p>0$ and require $\SA u = \SL$ in $\LqOYas,$ for some $q>0$ (i.e. use test functions in $\LqOY$). In this case, the proof of \cref{thm:12} would be nearly identical, as the space $\D$ of test functions used there 
is a subset of $\LqOY$ for all $q>0.$ One could also develop analogues of \cref{thm:3,thm:11,lem:svarwelldefined} in this setting---see e.g.~\cite[Theorem 3.20]{Gi:10} for an example of this approach for the stationary diffusion equation with lognormal diffusion coefficient.
\ere

\bre[Non-reliance on the Lax-Milgram theorem]
The above results hold for an arbitrary sesquilinear form and hence are applicable to a wide variety of PDEs; their main advantage is that they apply to PDEs whose stochastic variational formulations are not coercive. For example, as noted in \cref{sec:intro}, for the stationary diffusion equation \eqref{eq:diffusion} with coefficient $\kappa$ bounded uniformly below in $\omega,$ the bilinear form of \cref{prob:svar} is coercive; existence and uniqueness follow from the Lax-Milgram theorem, and hence the chain of results above leading to the well-posedness of \cref{prob:svar} is not necessary.
\ere

\bre[Overview of how these results are applied to the Helmholtz equation in \cref{sec:hhproof}]

We obtain the results for the Helmholtz equation via the following steps (which could also be applied to other SPDEs fitting into this framework):
\ben
\item\label[step]{it:step1pw} Define the map $\coeff$ (via $A,n,$ and $f$) such that for almost every $\omega \in \Omega$ there exists a solution of the deterministic Helmholtz EDP corresponding to $\coeff(\omega).$
\item\label[step]{it:step2pw} Define $u:\Omega\rightarrow X$ to map $\omega$ to the solution of the deterministic problem corresponding to $\coeff(\omega).$
\item Prove that \cref{con:coeffstoform,,con:A,,con:coeffstofunc,,con:L,,con:cborel,,con:C,,con:B,,con:K} hold, so that one can apply \cref{thm:3,,thm:11,,thm:12} along with \cref{lem:svarwelldefined,lem:uniq} to show \cref{prob:svsedp} is well-defined and $u$ is unique and satisfies \cref{prob:msedp,prob:somsedp,prob:svsedp}.
\een
\Cref{it:step1pw,it:step2pw} can be thought of as constructing a solution pathwise.
\ere

\section{Proof of the results in \cref{sec:general}}\label{sec:genproof}

A key ingredient in proving that the stochastic variational formulation is well-defined (\cref{lem:svarwelldefined}) is showing that the maps $\omega \mapsto \mleft[\mleft(\rhs\circ\coeff\mright)(\omega)\mright]\mleft(\vt(\omega)\mright)$ and $\omega \mapsto \mleft[\mleft[(\toform\circ\coeff(\omega)\mright]\mleft(\vo(\omega)\mright)\mright]\mleft(\vt(\omega)\mright)$ are measurable, for appropriate functions $\vo$ and $\vt$. Showing that these functions are measurable is not straightforward, because they both depend on $\omega$ in multiple places. However, the structure of the $\omega$-dependence in each case is similar, and so we first prove some general results that will be applicable to both of these cases.

\subsection{Preliminary lemmas}\label{sec:prelemmanary}
Throughout this \lcnamecref{sec:prelemmanary}, we assume we have two separable Banach spaces $\Zo$ and $\Zt$, and maps $\Pmap:\Omega\rightarrow \homspaceZoZt$ and $v:\Omega \rightarrow \Zo$. To simplify notation, we introduce the following definition.
\bde[Pairing map]
We define the map $\compPv:\Omega\rightarrow \Zt$ by
\beq\label{eq:comp}
\compPv(\omega) \de \mleft[\Pmap(\omega)\mright]\mleft(v(\omega)\mright).
\eeq
\ede

\bde[Product map]
Let $\toprodPv:\Omega\rightarrow \homspaceZoZt \times \Zo$ be defined by $\toprodPv(\omega) = \mleft(\Pmap(\omega),v(\omega)\mright).$
\ede

\ble[Product map is measurable]\label{lem:Pmeas}
When $\homspaceZoZt \times \Zo$ is equipped with the product topology, if $\Pmap$ and $v$ are measurable, then $\toprodPv:\Omega\rightarrow \homspaceZoZt\times \Zo$ is measurable.
\ele

\bpf[Proof of \cref{lem:Pmeas}]
By the result on the measurability of the Cartesian product of measureable functions (\cref{lem:measprod}), $\toprodPv$ is measurable with respect to $\big(\cF,\bigBorel{\homspaceZoZt}\otimes \Borel{\Zo}\big)$ (where $\cB$ denotes the Borel $\sigma$-algebra---see \cref{def:borelsigma}), as both of the coordinate functions $\Pmap$ and $v$ are measurable. Since $\homspaceZoZt$ and $\Zo$ are both metric spaces, they are both Hausdorff. As $\Zo$ is separable, \cref{lem:bogachev} on the product of Borel $\sigma$-algebras implies $\bigBorel{\homspaceZoZt} \otimes \Borel{\Zo} = \bigBorel{\homspaceZoZt \times \Zo}.$ Hence $\toprodPv$ is measurable with respect to $\big(\cF,\bigBorel{\homspaceZoZt\times \Zo}\big).$
\epf

\bde[Evaluation map]
The function $\evalZoZt:\homspaceZoZt \times \Zo \rightarrow \Zt$ is defined by
\beq\label{eq:evaldef}
\evalZoZt\big(\mleft(\cH,v\mright)\big) \de \cH(v) \quad\tfor \cH \in \homspaceZoZt \tand v \in \Zo.
\eeq
\ede

Observe that the pairing, product, and evaluation maps ($\compPv, \toprodPv,$ and $\evalZoZt$ respectively) are related by $\compPv = \evalZoZt \circ \toprodPv.$

\ble[Evaluation map is continuous]\label{lem:vstarcont}
The map $\evalZoZt$ is continuous with respect to the product topology on $\homspaceZoZt \times \Zo$ and the norm topology on $\Zt.$
\ele

The proof of \cref{lem:vstarcont} is straightforward and omitted.

\ble[Pairing map is measurable] \label{lem:abstractpairingborel}
If  $\Pmap$ and $v$ are measurable, then $\compPv$ is measurable.
\ele

\begin{proof}[Proof of \cref{lem:abstractpairingborel}]
By \cref{lem:Pmeas} $\toprodPv$ is measurable and by \cref{lem:vstarcont} $\evalZoZt$ is continuous. Therefore \cref{lem:contplusmeas} implies that $\compPv = \evalZoZt \circ \toprodPv$ is measurable.
 \end{proof}

\ble[$\mleft(\fullrhs\mright)(v)$ is measurable] \label{lem:thetaborel}
Under \cref{con:coeffstofunc,con:cborel}, for any measurable $\vt:\Omega\rightarrow Y,$ the function $\omega \mapsto \mleft[(\rhs\circ\coeff(\omega)\mright]\mleft(v(\omega)\mright)$ is measurable.
\ele

\begin{proof}[Proof of \cref{lem:thetaborel}]
The map $\coeff$ is measurable (by \cref{con:cborel}) and $\rhs$ is continuous  (by \cref{con:coeffstofunc}), therefore \cref{lem:contplusmeas} implies that $\fullrhs $ is measurable. Applying \cref{lem:abstractpairingborel} with $\Zo = Y,$ $\Zt=\FF$ (because $\Ys = \homspaceYF$), $\Pmap = \rhs \circ \coeff$, and $v = \vt,$ the result follows.
\end{proof}

\ble[$((\toform\circ\coeff)(\vo))(\vt)$
  is measurable] \label{lem:gammaborel}
If \cref{con:coeffstoform,con:cborel} hold and $\vo:\Omega \rightarrow X$ and $\vt:\Omega \rightarrow Y$ are measurable, then the function $\omega \mapsto \mleft[\mleft[\toform\circ\coeff(\omega)\mright]\mleft(\vo(\omega)\mright)\mright]\mleft(\vt(\omega)\mright)$ is measurable.
\ele

\begin{proof}[Proof of \cref{lem:gammaborel}]
  Since \cref{con:coeffstoform,con:cborel} hold, $\toform\circ\coeff$ is measurable by \cref{lem:contplusmeas}. Therefore by \cref{lem:abstractpairingborel} with $\Zo = X,$ $\Zt = \Ys,$ $\Pmap = \toform\circ\coeff$ and $v = \vo$, the map $\omega \rightarrow \mleft[\toform\circ\coeff(\omega)\mright]\mleft(\vo(\omega)\mright)$ is measurable. Therefore applying \cref{lem:abstractpairingborel} again with $\Zo = Y,$ $\Zt = \FF$, $\Pmap(\omega) = \mleft[\toform\circ\coeff(\omega)\mright]\mleft(\vo(\omega)\mright)$, and $v = \vt$, the result follows.
 \end{proof}

\subsection{Proofs of \cref{thm:3,thm:11,thm:12,lem:svarwelldefined,lem:uniq}}
\bpf[Proof of \cref{thm:3}]
We need to show $u:\Omega \rightarrow X$ is strongly measurable, satisfies the bound \eqref{eq:sbresult}, and therefore is Bochner integrable and is in the space $\LtOX.$ Our plan is to use \cref{cor:bochnersimple} to show $u$ is Bochner integrable, and establish \eqref{eq:sbresult} as a by-product. Since $u$ solves \cref{prob:meas}, $u$ is measurable. As $X$ is separable, it follows from \cref{cor:pettis} that $u$ is strongly measurable.
Define $N:X \rightarrow \RR$ by
$ N(v) \de \N{v}_X^2.$
Since $N$ is continuous, \cref{lem:contplusmeas} implies $N \circ u:\Omega \rightarrow \RR$ is measurable. 
Therefore, since both the left- and right-hand sides of \eqref{eq:sbe1} are measurable and \eqref{eq:sbe1} holds for almost every $\omega \in \Omega$ we can integrate \eqref{eq:sbe1} over $\Omega$ with respect to $\PP$ and obtain
\beq\label{eq:sbmid}
\int_\Omega \NX{u(\omega)}^2 \dd\PP(\omega) \leq \sum_{j=1}^m \NLoO{\Cj\fj},
\eeq
the right-hand side of which is finite since \cref{con:B} includes that $\Cj\fj \in \LoO$ for all $j = 1,\ldots,m.$ Since $u$ is strongly measurable, the bound \eqref{eq:sbmid} and \cref{cor:bochnersimple} with $p=2$ imply that $u$ is Bochner integrable. The norm $\NLtOX{u}$ is thus well-defined by \cref{def:bochnernorm} and \eqref{eq:sbmid} shows that \eqref{eq:sbresult} holds, and so in particular $\NLtOX{u} < \infty.$
\epf

\bpf[Proof of \cref{lem:svarwelldefined}]
We must show that for any $\vo \in \LtOX$ and any $\vt \in \LtOY$:
\bit
\item The quantities $\big[\Acomega \vo(\omega)\big]\big(\vt(\omega)\big)$ and $\Lcomega\big(\vt(\omega)\big)$ are Bochner integrable, so that the definitions of $\SA$ and $\SL$ as integrals over $\Omega$ make sense.
\item The maps $\SA(\vo)$ and $\SL$ are linear and bounded on $\LtOY,$ that is, $\SA:\LtOX\rightarrow\LtOYas$ and $\SL \in \LtOYas.$
\eit
It follows from these two points that $\SA$ and $\SL$ are well-defined.
Thanks to the groundwork laid in \cref{sec:prelemmanary}, $\big[\Acomega \vo(\omega)\big]\big(\vt(\omega)\big)$ and $\Lcomega\big(\vt(\omega)\big)$  are measurable by \cref{lem:gammaborel,lem:thetaborel} (which need \cref{con:coeffstoform,con:coeffstofunc,con:C}).
Their $\PP$-essential separability follows from \cref{con:coeffstoform,con:coeffstofunc,con:C,lem:esssep} and thus their strong measurability follows from \cref{cor:pettis} on the equivalence of measurability and strong measurability when the image is separable. Their Bochner integrability then follows from the Bochner integrability condition in \cref{thm:bochnercond} (with $V=\FF$) and the Cauchy--Schwarz inequality since
\begin{align}
\int_\Omega \abs{\Lcomega\big(\vt(\omega)\big)}\dd\PP(\omega)&\leq \int_\Omega \NYs{\mleft(\rhs\circ\coeff\mright)(\omega)}\NY{\vt(\omega)}\dd\PP(\omega)\nonumber\\
&\leq \NLtOYs{\rhs\circ\coeff}\NLtOY{\vt},\label{eq:Lfirst}
\end{align}
which is finite by \cref{con:L}, and 
\begin{align}
\int_\Omega \Big|\big[\Acomega \vo(\omega)\big]\big(\vt(\omega)\big)\Big|\dd\PP(\omega) 
&\leq \esssup_{\omega \in \Omega} \Nhomspace{\Acomega} \int_\Omega \NX{\vo(\omega)}\NY{\vt(\omega)}\dd\PP(\omega)\nonumber\\
&\leq \NLiOhomspace{\toform\circ\coeff}\NLtOX{\vo}\NLtOY{\vt},\label{eq:Afirst}
\end{align}
which is finite by \cref{con:A}.

We now show $\SL\in  \LtOYas$ and $\SA:\LtOX \rightarrow \LtOYas.$ Observe that 

\noindent $\abs{\SL(\vt)} \leq \int_\Omega \abs{\Lcomega\mleft(\vt(\omega)\mright)}\dd\PP(\omega)$ and $\abs{\mleft[\SA \mleft(\vo\mright)\mright](\vt)} \leq \int_\Omega \abs{\mleft[\Acomega \vo(\omega)\mright]\mleft(\vt(\omega)\mright)}\dd\PP(\omega)$ and thus by \eqref{eq:Lfirst} and \eqref{eq:Afirst} $\SL$ and $\SA(\vo)$ are bounded. They are clearly linear, and so it follows that $\SL \in \LtOYas$ and $\SA(\vo)\in \LtOYas,$ i.e., $\SA:\LtOX \rightarrow \LtOYas.$
\epf

\bpf[Proof of \cref{thm:11}]
In order to show that $u$ solves \cref{prob:svar}, we must show:
\ben
\item\label[point]{it:111} either the functional $\SLc \in \LtOYas$ or the functional $\SA(u) \in \LtOYas$, and
\item\label[point]{it:112} the equality \eqref{eq:stoeq} holds.
\een

For \cref{it:111} we show that $\SL \in \LtOYas,$ (since this is easier than showing $\SA(u) \in \LtOYas$); in fact the proof of this is contained in the proof of \cref{lem:svarwelldefined}.

For \cref{it:112}, since $u$ solves \cref{prob:lt}, for $\PP$-almost every $\omega \in \Omega$ we have
$\Acomega u(\omega) = \Lcomega$
in $\Ys.$ Hence, for any $v \in \LtOY$ we have
\beq\label{eq:midwaytoeu2}
\big[\Acomega u(\omega)\big]\big(v(\omega)\big) = \Lcomega\big(v(\omega)\big)
\eeq
for $\PP$-almost every $\omega \in \Omega.$ Since $\SL \in \LtOYas$, the right-hand side of \eqref{eq:midwaytoeu2} is a strongly measurable function with finite integral. Hence the left-hand side of \eqref{eq:midwaytoeu2} is as well, and we can integrate over $\Omega$ to conclude
$\big[\SAc u\big](v) = \SLc(v) \tforall v \in \LtOY,$
that is, $\SAc u = \SLc$ in $\LtOYas.$
\epf
The following lemma is needed for the proof of \cref{thm:12}.
\ble\label{lem:settheory}
Let $\diff:\Omega\times Y \rightarrow \FF.$  For $y \in Y,$ define $\Omegay \de \set{\omega \in \Omega \st \diff(\omega,y)=0}$ and define $\Omegat \de \set{\omega \in \Omega \st \diff(\omega,y)=0 \tforall y \in Y}.$ If
\bit
\item for all $\omega \in \Omega,$ $\diff(\omega,\cdot)$ is a continuous functional on $Y$ and
\item for all $y \in Y,$ the map $\diff(\cdot,y):\Omega\rightarrow \FF$ is measurable and $\PP(\Omegay)=1,$
  \eit
  then $\PP(\Omegat)=1.$
\ele

\bpf[Proof of \cref{lem:settheory}]
We must show that the set $\Omegat \in \cF,$ and $\PP(\Omegat)=1.$ Observe that, for any $y \in Y$, the set $\Omegay \in \cF,$ since  $\Omegay = \delta(\cdot,y)^{-1}\mleft(\set{0}\mright),$ which  is the preimage under a measurable map of a measurable set. 

Since $Y$ is a Hilbert space, it is separable, and therefore it has a countable dense subset $\mleft(\yn\mright)_{n \in \NN}.$ We will show that $\PP\mleft(\cap_{n \in \NN} \Omegayn\mright)=1$ and $\Omegat = \cap_{n \in \NN} \Omegayn.$ The set $\cap_{n \in \NN} \Omegayn \in \cF,$ as $\cF$ is a $\sigma$-algebra and $\PP\mleft(\cup_{n \in \NN} \Omegaync\mright) \leq \sum_{n \in \NN} \PP\mleft(\Omegaync\mright) = 0,$ and hence $\PP\mleft(\cap_{n\in\NN} \Omegayn\mright)=1.$ To next show $\Omegat =  \cap_{n \in \NN} \Omegayn$ we observe that $\Omegat = \cap_{y \in Y} \Omegay$ and $\cap_{y \in Y} \Omegay \subseteq \cap_{n \in \NN} \Omegayn.$ It therefore suffices to show $\cap_{n \in \NN} \Omegayn \subseteq \cap_{y \in Y} \Omegay$ to conclude $\Omegat =  \cap_{n \in \NN} \Omegayn.$

Fix $y \in Y.$ By density of $\mleft(\yn\mright)_{n \in \NN}$, there exists a subsequence $\ynmseq$ such that $\ynm \rightarrow y$ as $m \rightarrow \infty.$ Fix $\omega \in \cap_{n \in \NN} \Omegayn.$ Note that $\omega \in \cap_{m \in \NN} \Omegaynm;$ that is, for all $m \in \NN,$ $\diff(\omega,\ynm) =0.$ As $\diff(\omega,\cdot)$ is a continuous function on $Y$, $\diff(\omega,\ynm) \rightarrow \diff(\omega,y)$ as $m \rightarrow \infty.$ But as previously noted, $\diff(\omega,\ynm)=0$ for all $m \in \NN.$ Hence we must have $\diff(\omega,y)=0,$ and thus $\omega \in \Omegay.$ Since $\omega \in \cap_{n \in \NN} \Omegayn$ was arbitrary, it follows that $\cap_{n \in \NN} \Omegayn \subseteq \Omegay,$ and since $y \in Y$ was arbitrary, it follows that $\cap_{n \in \NN} \Omegayn \subseteq \cap_{y \in Y} \Omegay$ as required.
\epf

\bpf[Proof of \cref{thm:12}]
Let $u \in \LtOX$ solve \cref{prob:svar}. We need to show that $u$ solves \cref{prob:lt}. Observe that $u$ solving \cref{prob:lt} means $\Acomega(u(\omega)) = \mleft(\Lcomega\mright)(\omega)$ in $\Ys$ for almost every $\omega \in \Omega.$ We now use an idea from \cite[Theorem 3.3]{Gi:10}. Our plan is to use test functions of the form $y\Ind{E},$ where $y \in Y$ and $E \in \cF$ to reduce \cref{prob:svar} to the statement
\beqs
\int_E \mleft[\Acomega\big(u(\omega)\big)\mright]\big(y(\omega)\big) \dd\PP(\omega) = \int_E \mleft[\mleft(\Lcomega\mright)(\omega)\mright]\big(y(\omega)\big) \dd\PP(\omega)\quad \tforall E \in \cF
\eeqs
and then show this implies $u$ satisfies \cref{prob:lt} via \cref{lem:gotoae}.

First define the space
$\D := \set{y\Ind{E} \st y \in Y, E \in \cF}.$
It is straightforward to see that the elements of $\D$ are maps from $\Omega$ to $Y.$ The fact that $\D \subseteq \LtOY$ follows via the following three steps:

\ben
\item The elements of $\D$ are measurable, indeed the indicator function of a measurable set is a measurable function $\Omega\rightarrow\RR,$ and multiplication by $y \in Y$ is a continuous function $\RR\rightarrow Y.$ Hence elements of $\D$ are measurable by \cref{lem:contplusmeas}.
\item As $Y$ is a separable Hilbert space, it follows from \cref{cor:pettis} that the elements of $\D$ are strongly measurable.
\item $\NLtOY{y\Ind{E}} = \sqrt{\PP\mleft(E\mright)}\NY{y} < \infty$ for all $y \in Y, E \in \cF.$
  \een

Since \cref{prob:svar} is well-defined, and $u$ solves \cref{prob:svar}, and $\D \subseteq \LtOY,$ we have that $\mleft[\SA u\mright](v) = \SL(v) \tforall v \in \D.$ Therefore, we have

\beq\label{eq:initialint}
\int_\Omega \mleft[\Acomega\mleft(u(\omega)\mright)\mright]\mleft(y\Ind{E}(\omega)\mright) \dd\PP(\omega) = \int_\Omega \mleft[\Lcomega\mright]\mleft(y\Ind{E}(\omega)\mright) \dd\PP(\omega)
\eeq
for all $y \in Y$ and $E \in \cF.$ If we define $\diff:\Omega\times Y \rightarrow \FF$ by $\diff(\omega,y) \de \mleft[\Acomega\mleft(u(\omega)\mright) - \Lcomega\mright]\mleft(y\mright)$ then, by the definition of $\Ind{E},$ \eqref{eq:initialint} becomes
\beq\label{eq:intoverE}
\int_E \diff(\omega,y) \dd\PP(\omega)=0\quad \tforall E \in \cF.
\eeq
To conclude $u$ solves \cref{prob:lt} we must show $\diff(\omega,y)=0$ for all $y \in Y,$ almost surely. We will use \cref{lem:gotoae}, so the first step is to show that for all $y \in Y$ $\diff(\cdot,y)$ is Bochner integrable. This follows from the fact that \cref{prob:svar} is well-defined, and thus the quantities $\big[\Acomega \vo(\omega)\big]\big(\vt(\omega)\big)$ and $\Lcomega\big(\vt(\omega)\big)$ are Bochner integrable for any $\vo\in \LtOX,\vt \in \LtOY.$ In particular, they are Bochner integrable when $\vo=u,$ and $\vt=y\Ind{E}$ and thus their difference $\diff$ is Bochner integrable. Secondly, $\diff(\omega,\cdot)$ is a continuous function on $Y$ since $\Acomega\mleft(u(\omega)\mright)$, $\mleft(\Lcomega\mright)(\omega) \in \Ys,$ for all $\omega \in \Omega.$

We now show $\diff(\omega,y)=0$ for all $y \in Y,$ almost surely. For $y \in Y$ define the set $\Omegay \de \set{\omega \in \Omega \st \diff(\omega,y)=0};$ by \eqref{eq:intoverE} and \cref{lem:gotoae} we have that $\PP(\Omegay)=1$ for all $y \in Y.$ By \cref{lem:settheory}, $  \diff(\omega,y)=0$ for all $y \in Y$, almost surely, that is, $\Acomega u(\omega) = \Lcomega$ almost surely; it follows that $u$ solves \cref{prob:lt}.
\epf

\bre[Connection with the argument in {\cite[Remark 2.2]{MuSt:11}}]
The argument in 

\noindent \cref{lem:settheory} and the final part of \cref{thm:12} closely mirrors the result in \cite[Remark 2.2]{MuSt:11}. Indeed, we prove in general that
\beqs
\PP\big(\diff(\omega,y)=0\big)=1\text{ for all } y \in Y \quad\text{implies} \quad\PP\big(\diff(\omega,y)=0\text{ for all } y \in Y\big)=1,
\eeqs
and \cite[Remark 2.2]{MuSt:11} shows an analogous result for the stationary diffusion equation \eqref{eq:diffusion} with non-uniformly coercive and unbounded coefficient $\kappa.$
\ere

\bpf[Proof of \cref{lem:uniq}]
\emph{Proof of \cref{part:uniq1}.} Suppose $\uo,\ut:\Omega\rightarrow X$ solve \cref{prob:meas}. Let $E = \set{\omega \in \Omega \st \uo(\omega) \neq \ut(\omega)}.$ Denote by $\Eo$ and $\Et$ the sets (of measure zero) where the variational problems for $\uo$ and $\ut$ fail to hold, i.e. $\Eo,\Et \in \cF$ with $\PP(\Eo)=\PP(\Et)=0$ and 
\beqs
\Acomega\mleft(\uo(\omega)\mright) \neq \Lcomega \tiff \omega \in \Eo,\quad\text{and}\quad \Acomega\mleft(\ut(\omega)\mright) \neq \Lcomega \tiff \omega \in \Et.
\eeqs As  $\ker\mleft(\Acomega\mright) = \set{0}$ $\PP$-almost surely, there exists $\Eth \in \cF$ such that $\PP(\Eth) = 0$ and

\noindent $\ker\mleft(\Acomega\mright) \neq \set{0} \tiff \omega \in \Eth.$
We claim $E \subseteq \Eo\cup\Et\cup\Eth.$ Indeed, if $\uo(\omega) \neq \ut(\omega)$ then either: (i) at least one of $\uo$ and $\ut$ does not solve \cref{prob:meas} at $\omega$ or (ii) $\uo$ and $\ut$ both solve \cref{prob:meas} at $\omega,$ but $\ker\mleft(\Acomega\mright) \neq \set{0}.$
Since $\PP(E_j)=0, j = 1,2,3,$ we have $\PP(\Eo\cup\Et\cup\Eth) = 0.$  Therefore $E \in \cF$ and $\PP(E)=0$ since $\OFP$ is a complete probability space; hence $\uo = \ut$ almost surely, as required.

\emph{Proof of \cref{part:uniq2}.} By \cref{rem:imm}, if $\uo,\ut \in \LtOX$ solve \cref{prob:lt}, then all the representatives of the equivalence classes of $\uo$ and $\ut$ solve \cref{prob:meas}. Hence, by \cref{part:uniq1}, any representative of $\uo$ and any representative of $\ut$ differ only on some set (depending on the representatives) of $\PP$-measure zero in $\Omega.$ Therefore $\uo=\ut$ in $\LtOX,$ by definition of $\LtOX.$

\emph{Proof of \cref{part:uniq3}.} As \cref{prob:svar} is well-defined, by \cref{rem:imm,thm:12}, if $\uo$ and $\ut$ solve \cref{prob:svar}, then $\uo$ and $\ut$ also solve \cref{prob:meas}. We then repeat the reasoning in the proof of \cref{part:uniq2} to show $\uo=\ut$ in $\LtOX.$
\epf

\section{Proofs of \cref{thm:hh-gen,thm:hh-hetero}}\label{sec:hhproof}
In \cref{sec:placing} we place the Helmholtz stochastic EDP into the framework developed in \cref{sec:gen-framework}. In \cref{sec:hh-cond} we give sufficient conditions for the Helmholtz stochastic EDP to satisfy \cref{con:coeffstoform,,con:coeffstofunc,,con:cborel}, etc.. In \cref{sec:applying} we apply the general theory developed in \cref{sec:gen-framework} to prove \cref{thm:hh-gen,,thm:hh-hetero}.

\subsection{Placing the Helmholtz stochastic EDP into the framework of \cref{sec:gen-framework}}\label{sec:placing}
Recall $R>0$ is fixed. We let $X=Y=\HozDDR$ and define the norm $\NW{v}^2 \de \NLtDR{\grad v}^2 + k^2 \NLtDR{v}^2$ on $\HozDDR.$ Throughout this section, $\Az,\nz,$ and $\fz$ will be deterministic functions. Recall that since the supports of $1-n,$ $I-A,$ and $f$ are compactly contained in $\BR$, we can consider $A, n,$ and $f$ as functions on $\DR$ rather than on $\Dp$. In order to define the space $\cC$ and the maps $\coeff,\toform,$ and $\rhs$ we define the following function spaces on $\DR$.

\bde[Compact-support spaces]\label{def:compsuppspace}
Let
\beqs
\LtRDR \de \set{\fz \in \LtDR \st \supp\mleft(\fz\mright) \compcont \BR}.
\eeqs
\begin{align*}
\LiRminDRRR \de \big\{&\nz \in \LiDRRR \st\supp\mleft(1-\nz\mright) \compcont \BR,\\
&\text{there exists } \alphanz > 0 \text{ such that } \nz(\bx) \geq \alphanz \text{ almost everywhere }\big\},\\
\hspace{-3cm}\LiRminDRRRdtd \de \Big\{&\Az \in \LiDRRRdtd \st\Az(\bx) \text{ is symmetric almost everywhere,}\\
&\supp\mleft(I-\Az\mright) \compcont \BR, \text{there exists } \alphaAz>0 \text{ s.~t.~}\alphaAz \leq \Az(\bx)\\
&\text{almost everywhere, in the sense of quadratic forms}\Big\},\text{ and}\\
&\hspace{-3cm}\WoiRminDRRRdtd \de \set{\Az \in \LiRminDRRRdtd \st\Az \in \WoiDRRRdtd}.
\end{align*}
\ede

Observe that the norm on $\LiDRRR$ induces a metric on $\LiRminDRRR,$ and similarly for $\LiRDRRRdtd,$ $\WoiRminDRRRdtd,$ and $\LtRDR.$ These spaces are not vector spaces, and are not complete, but completeness and being a vector space are not required in what follows---we only need them to be metric spaces.

\bde[Deterministic form and functional]

\noindent For $\mleft(\Az,\nz,\fz\mright) \in \LiRDRRRdtd \times \LiRminDRRR \times \LtRDR$ let the sesquilinear form $\aGm$ on $\HozDDR \times \HozDDR$  and the antilinear functional $\Lh$ on $\HozDDR$ be given by
\begin{align*}
\aGm\mleft(\vo,\vt\mright) &\de \int_{D_R} \Big(\mleft(\Az \grad \vo\mright)\cdot \grad \vtb \rangle 
 - k^2 \nz\, \vo\,\vtb \Big)\dd\Leb- \big\langle T_R \gamma \vo,\gamma \vt\big\rangle_{\Gamma_R}, \quad\text{and}\\
 \Lh(\vt) &\de \int_{D_R} \fz\, \vtb\,\dd\Leb, \quad\text{ for } \vo, \vt \in \HozDDR.
\end{align*}
\ede

\bprob[Helmholtz EDP]\label{prob:edp}
For $\mleft(\Az,\nz,\fz\mright) \in \LiRDRRRdtd \times \LiRDRRR \times \LtRDR$ find 
$\uz \in \HozDDR$ such that $\aGm(\uz,v) = \Lh(v) \tforall v \in \HozDDR.$
\eprob

\bde[$\dinfty$ metric]\label{def:dinfty}
Let $\mleft(\Xo,\done\mright),\ldots,\mleft(\Xm,\dm\mright)$ be metric spaces. The \defn{$\dinfty$ metric} on the Cartesian product $\Xo\times\cdots\times\Xm$ is defined by
\beqs
\dinfty\mleft(\mleft(\xo,\ldots,\xm\mright),\mleft(\yo,\ldots,\ym\mright)\mright) \de \max_{j = 1,\ldots,m} \dmetj\mleft(\xj,\yj\mright).
\eeqs
\ede

\bde[The input space $\cC$]\label{def:cCHh}
We let
$\cC \de \WoiRminDRRRdtd \times \LiRminDRRR \times \LtRDR$
with topology given by the $\dinfty$ metric.
\ede

\bde[The input map $\coeff$]\label{def:inputmap}
Define $\coeff:\Omega\rightarrow\cC$ by
$\coeff(\omega) = \mleft(A(\omega),n(\omega),f(\omega)\mright).$
\ede

\bde[The maps $\toform$ and $\rhs$ for the Helmholtz stochastic EDP]
Let
\beq\label{eq:hhform}
\toform\mleft(\mleft(\Az,\nz,\fz\mright)\mright) \de \aGm \quad \text{and} \quad \rhs\mleft(\mleft(\Az,\nz,\fz\mright)\mright) \de \Lh,
\eeq
where the definition of $\toform$ is understood in terms of the equivalence between $\homspace$ and sesquilinear forms on $X \times Y.$
\ede

\subsection{Verifying the Helmholtz stochastic EDP satisfies the general conditions in \cref{sec:gen-framework}}\label{sec:hh-cond}

\ble[{\Cref{con:cborel,con:C}} for Helmholtz stochastic EDP]\label{lem:hh-borelC}
If $A,n,$ and $f$ are strongly measurable, then $\coeff$ defined by  \cref{def:inputmap} satisfies \cref{con:cborel,con:C}.
\ele

\bpf
Since $A,n,$ and $f$ are strongly measurable, by \cref{thm:pettis} they are measurable and $\PP$-essentially separably valued. By \cref{lem:measprod}, it follows that $\coeff$ is measurable, so $\coeff$ satisfies \cref{con:cborel}. By \cref{lem:prodsep}, it follows that $\coeff$ is $\PP$-essentially separably valued, so $\coeff$ satisfies \cref{con:C}.
\epf

\ble[\Cref{con:coeffstoform,con:coeffstofunc} for Helmholtz stochastic EDP]\label{lem:hh-AL}
The maps $\toform$ and $\rhs$ given by \eqref{eq:hhform} satisfy \cref{con:coeffstoform,con:coeffstofunc}.
\ele

\bpf[Proof of \cref{lem:hh-AL}]
We need to show that if $(\Am,\nm,\fm) \rightarrow (\Az,\nz,\fz)$ in $\cC$ then $\toform((\Am,\nm,\fm)) \rightarrow \toform((\Az,\nz,\fz))$ in $\homspace,$ and similarly for $\rhs.$ We have, for $\vo \in X,\vt \in Y,$
\begin{align*}
&\biggabs{\Big[\big[\toform\mleft(\Am,\nm,\fm\mright) - \toform\mleft(\Az,\nz,\fz\mright)\big]\mleft(\vo\mright)\Big]\mleft(\vt\mright)}\\
&\quad\quad\quad\,\,= \abs{\int_{\DR} \Big(\big(\mleft(\Am-\Az\mright)\grad \vo\big)\cdot\grad\vtb- k^2 \mleft(\nm-\nz\mright) \vo\vtb\Big)\dd\Leb}\\
&\quad\quad\quad\,\,\leq\NLiDRRRdtd{\Am-\Az}\NLtDR{\grad \vo}\NLtDR{\grad \vt} 
\\
&\quad\quad\quad\,\,\quad\quad
+k^2\NLiDRRR{\nm-\nz}\NLtDR{\vo}\NLtDR{\vt}\\
&\quad\quad\quad\,\,\leq 2\dinfty((\Am,\nm,\fm),(\Az,\nz,\fz))\Nw{\vo}\Nw{\vt},
\end{align*}
Hence if $(\Am,\nm,\fm) \rightarrow (\Az,\nz,\fz)$ in $\cC,$ then $\toform((\Am,\nm,\fm)) \rightarrow \toform((\Az,\nz,\fz))$ in 

\noindent $\homspace.$ We also have
\beqs
\Bigabs{\big[\rhs\mleft(\mleft(\Am,\nm,\fm\mright),\mright) - \rhs\mleft(\mleft(\Az,\nz,\fz\mright)\mright)\big]\mleft(\vt\mright)} 
 = \abs{\int_{\DR} \mleft(\fm- \fz\mright)\vtb\,\dd\Leb}
 \leq \NLtDR{\fm-\fz} \frac{\Nw{\vt}}{k}.
\eeqs
Hence if $(\Am,\nm,\fm) \rightarrow (\Az,\nz,\fz)$ in $\cC,$ then $\rhs((\Am,\nm,\fm)) \rightarrow \rhs((\Az,\nz,\fz))$ in $\Ys.$
\epf

\bde[The solution operator $\sol$]
Define $\sol:\cC\rightarrow \HozDDR$ by letting 

\noindent $\sol\mleft(\Az,\nz,\fz\mright) \in \HozDDR$ be the solution of the Helmholtz EDP (\cref{prob:edp}).
\ede

\bth[$\sol$ is well defined]\label{thm:deteu}
For $\mleft(\Az,\nz,\fz\mright) \in \cC$ the solution $\sol\mleft(\mleft(\Az,\nz,\fz\mright)\mright)$ of the Helmholtz EDP (\cref{prob:edp}) exists, is unique, and depends continuously on $\fz.$
\enth

\bpf[Proof of \cref{thm:deteu}] Since $\Real{-\langle \TrR \gamma v,\gamma v\rangle_{\GR}} \geq 0$ for all $v \in \HozDDR$ (see, e.g.~\cite[Theorem 2.6.4]{Ne:01}), $\aGm$ satisfies a G\r{a}rding inequality. Since the inclusion $\HozDDR \hookrightarrow \LtDR$ is compact, Fredholm theory shows that uniqueness implies well-posedness (see, e.g.~\cite[Theorem 2.34]{Mc:00}). Since $A$ is Lipschitz and $n$ is $L^\infty$, uniqueness follows from the unique continuation results in  \cite{JeKe:85,GaLi:87};
see \cite[Section 2]{GrSa:18} for these results specifically applied to Helmholtz problems.
\epf

\ble[Continuity of solution operator for Helmholtz stochastic EDP]\label{lem:solcont}
For the 

\noindent  Helmholtz stochastic EDP, the solution operator $\sol:\cC\rightarrow\HozDDR$ is continuous.
\ele

\bpf[Sketch Proof of \cref{lem:solcont}]
Let $(\Az,\nz,\fz), (\Ao,\no,\fo) \in \cC,$ with $\sol((\Az,\nz,\fz)) = \uz$ and $\sol((\Ao,\no,\fo)) = \uo.$
Then for any $v \in \HozDDR$ we have, for $j=0,1,$

\beqs\label{eq:S1}
\mleft[\mleft[\toform((A_j, n_j, f_j))\mright](u_j)\mright](v) = \mleft[\rhs((A_j,n_j,f_j))\mright](v).
\eeqs
Continuity of $\sol$ then follows from:

\ben
\item Deriving the Helmholtz equation with coefficients $\Az$ and $\nz$ satisfied by $\ud\de\uz-\uo.$
\item\label[point]{it:contrecall} Recalling that the well-posedness result of \cref{thm:deteu} holds when $\fz \in \LtRDR$ is replaced by a right-hand side in $(\HozDDR)^*$; see, e.g., \cite[Theorem 2.34]{Mc:00}.
\item Applying the result in \cref{it:contrecall} to obtain a bound
  $\NW{\ud} \leq C(\Az,\nz) \NHozDDRas{F}.$

\item Showing $\NHozDDRas{F}$ depends on $\NLtDR{\grad \uo},$ $\NLtDR{\uo},$ $\NLiDRRRdtd{\Ao-\Az},$ $\NLiDRRR{\no-\nz},$ and $\NLtD{\fz-\fo}.$
\item Eliminating the dependence on $\uo$ by writing $\uo = \uz-\ud$ and moving terms in $\ud$ to the left-hand side, to obtain a bound on $\ud$ of the form
  \begin{align*}
&\NLtDR{\grad \ud} + k\NLtDR{\ud} \\
&\qquad\leq \Ctilde\Big(\uz,\Az,\nz,\NLiDRRRdtd{\Ao-\Az},\NLiDRRR{\no-\nz},\NLtDR{\fz-\fo}\Big).
\end{align*}
\item Concluding that $\ud\rightarrow0$ in $\HozDDR$ as $(\Ao,\no,\fo) \rightarrow (\Az,\nz,\fz)$ in $\cC.$
\een\epf

\ble[\Cref{con:K} for the Helmholtz stochastic EDP]\label{lem:hh-K}
The Helmholtz stochastic EDP satisfies \cref{con:K}.
\ele

\bpf[Proof of \cref{lem:hh-K}]
This condition holds immediately from \cref{thm:deteu}.
\epf

To prove that \cref{con:B} holds for the Helmholtz stochastic EDP, we first state the deterministic analogues of \cref{con:hh-hetero,thm:hh-hetero}.

\bcon[Nontrapping condition for Helmholtz EDP {\cite[Condition 2.4]{GrPeSp:18}}]\label{cond:1}
$d=2,3$, $\Dm$ is star-shaped with respect to the origin, $\Az\in \WoiDRRRdtd$, $\nz\in \WoiDRRR$, and there exist $\consto, \constt>0$ such that,
for almost every $\bx\in \Dp$,
\begin{align}\label{eq:A1}
&\Az(\bx) - (\bx\cdot\nabla)\Az(\bx) \geq \consto, \text{ in the sense of quadratic forms, and }\\
\label{eq:n1}
& \hspace{2cm}\nz(\bx)+ \bx\cdot\nabla \nz(\bx) \geq \constt.
\end{align}
\econ

\bth[Well-posedness of the Helmholtz EDP under \cref{cond:1} {\cite[Theorem 2.5]{GrPeSp:18}}]\label{thm:eubedp}
Let $(\Az,\nz,\fz) \in \cC$ and suppose $\Az$ and $\nz$ satisfy \cref{cond:1}. Then the solution of the Helmholtz EDP (\cref{prob:edp}) exists and is unique. Furthermore, given $\kz > 0$ for all $k \geq \kz,$ the solution $\uz$ of the Helmholtz EDP satisfies the bound 
\beq\label{eq:heterobound1}
\consto \NLtDR{\grad \uz}^2 + \constt k^2 \NLtDR{\uz}^2\leq \Co \NLtDR{\fz}^2,\,\, \text{where} \,\, \Co := 4\mleft[\frac{R^2}{\consto} + \frac{1}{\constt}\mleft(R+ \frac{d-1}{2\kz}\mright)^2\mright].
\eeq
\enth

We can now prove \cref{con:B} holds for the Helmholtz stochastic EDP.

\ble[\Cref{con:B} for Helmholtz stochastic EDP]\label{lem:hh-B}
If \cref{con:hh-fAn,con:hh-hetero} hold, then \cref{con:B} holds for the Helmholtz stochastic EDP.
\ele

\bpf[Proof of \cref{lem:hh-B}]
As \cref{con:hh-hetero} holds,  \cref{cond:1} holds for $\PP$-almost every $\omega \in \Omega$ (with $\Az = A(\omega),$ $\nz = n(\omega),$ $\consto = \muo(\omega),$ and $\constt = \mut(\omega)$). Hence, by \cref{thm:eubedp} the bound \eqref{eq:sbe1} holds for all $k \geq \kz$, with $X = \HozDDR, m=1,$
\beqs
\Co(\omega) = \frac4{\min\set{\muo(\omega),\mut(\omega)}}\mleft[\frac{R^2}{\muo(\omega)} + \frac{1}{\mut(\omega)}\mleft(R+ \frac{d-1}{2\kz}\mright)^2\mright],
\eeqs
and $\fo = \NLtDR{f(\omega)}^2.$ It now remains to show that $\Co\,\NLtDR{f}^2 \in \LoO.$ We first show $\Co\,\NLtDR{f}^2$ is measurable and then show that it lies in $\LoO.$ To show measurability, we rewrite $\Co(\omega)$ as 
\beqs
\Co(\omega) =  \max\left\{\frac{2 R^2}{\mu_1^2(\omega)} + \frac{2}{\mu_1(\omega)\mu_2(\omega)}\left(R+ \frac{d-1}{2\kz}\right)^2,\frac{2 R^2}{\mu_1(\omega)\mu_2(\omega)} + \frac{2}{\mu_2^2(\omega)}\left(R+ \frac{d-1}{2\kz}\right)^2\right\}.
\eeqs
The functions $\muo^{-1}$ and $\mut^{-1}$ are measurable by assumption; to conclude $\Co$ is measurable we use the facts (see e.g. \cite[Theorems 19.C, 20.A]{Ha:74}): (i) the square of a measurable function is measurable, and (ii) the product, sum, and maximum of two measurable functions are measurable. Under 
\cref{con:hh-fAn}, the function $f$ lies in the Bochner space $\LtOLtDR.$ Therefore, $f$ is strongly measurable and hence $f$ is measurable  by \cref{thm:pettis}. The map $f\mapsto\NLtDR{f}^2$ is clearly continuous, and therefore $\fo$ is measurable by \cref{lem:contplusmeas}. As the product of two measurable functions is measurable, it follows that $\Co\,\NLtDR{f}^2$ is measurable.

We now show that $\Co\NLtDR{f}^2 \in \LoO.$ The assumptions $1/\muo,1/\mut \in \LtO$ and the Cauchy--Schwarz inequality imply $1/(\muo\mut) \in \LoO.$ Therefore the maps,
\beqs\omega \mapsto \frac{2 R^2}{\mu_1^2(\omega)} + \frac{2}{\mu_1(\omega)\mu_2(\omega)}\left(R+ \frac{d-1}{2\kz}\right)^2 \text{ and } \omega \mapsto \frac{2 R^2}{\mu_1(\omega)\mu_2(\omega)} + \frac{2}{\mu_2^2(\omega)}\left(R+ \frac{d-1}{2\kz}\right)^2
\eeqs
are in $\LoO.$ Since the maximum of two functions in $\LoO$ is also in $\LoO,$ it follows that $\Co \in \LoO.$ \cref{con:hh-fAn} implies that $\NLtDR{f}^2 \in \LoO.$

To conclude $\Co\NLtDR{f}^2 \in \LoO,$ observe that the only dependence of $\Co$ on $\omega$ is through $\muo$ and $\mut.$ As $\muo$ and $\mut$ are assumed independent of $f,$ and measurable functions of independent random variables are independent \cite[p.236]{Lo:77} it follows that $\Co$ and $\NLtDR{f}^2$ are independent, and therefore
\begin{align}\nonumber
&\hspace{-2ex}\NLoO{\Co\NLtDR{f}^2} = \int_\Omega \Co(\omega) \NLtDR{f(\omega)}^2 \dd\PP(\omega) 
\\ &= \mleft(\int_\Omega \Co(\omega) \dd\PP(\omega)\mright)\mleft(\int_\Omega \NLtDR{f(\omega)}^2 \dd\PP(\omega)\mright)= \NLoO{\Co}\NLoO{\NLtDR{f}^2} < \infty.\label{eq:normsplit}
\end{align}
Therefore $\Co\NLtD{f}^2 \in \LoO$ as required. We take the expectation (equivalently, the $L^1$ norm) of \eqref{eq:heterobound1} (with $\Az = A(\omega)$ etc.) and use \eqref{eq:normsplit} to obtain \eqref{eq:Sbound1}.
\epf

\bre[The case when $f,$ $\muo,$ and $\mut$ are not independent]\label{rem:notindep}
\Cref{rem:planewave} shows that for the physically relevant example of scattering by a plane wave, $f,$ $\muo,$ and $\mut$ may not be independent. In this case, if we replace the requirements in \cref{con:hh-hetero} that $f \in \LtOLtD$ and $1/\muo,\,1/\mut \in \LtO$ with the stronger requirements $f \in \LfOLtD$ and $1/\muo,\,1/\mut \in \LfO$, then one can obtain the bound
\beqs
\NLtOHozDDR{\grad u}^2 + k^2\NLtOHozDDR{ u}^2\leq \NLtO{\Co} \NLfOLtDR{f}^2.
\eeqs
Indeed, instead of independence, we use the Cauchy--Schwartz inequality in \eqref{eq:normsplit} to conclude
\beqs
\NLoO{\Co \NLtDR{f}^2}  \leq \NLtO{\Co}\NLtO{\NLtDR{f}^2} = \NLtO{\Co} \NLfOLtDR{f}^2.
\eeqs
\ere

\ble[\Cref{con:L} for Helmholtz stochastic EDP]\label{lem:hh-L}
If $f \in \LtOLtDR$ and $A$ and $n$ are strongly measurable, then \cref{con:L} holds for the Helmholtz stochastic EDP.
\ele

\bpf[Proof of \cref{lem:hh-L}]
Since $A,n,$ and $f$ are strongly measurable, \cref{con:cborel,con:C} hold by \cref{lem:hh-borelC}; i.e., $\coeff$ is both measurable and $\PP$-essentially separably valued. Furthermore, by \cref{thm:pettis} $\coeff$ is strongly measurable. By \cref{lem:hh-AL}, \cref{con:coeffstofunc} holds, so the map $\rhs$ is continuous. Hence, by \cref{lem:contplusstrong}, $\rhs\circ\coeff$ is strongly measurable. We also have that
$\NYs{\mleft(\rhs\circ\coeff\mright)(\omega)} = \NLtDR{f(\omega)}/k,$
and thus $\rhs\circ\coeff \in \LtOYs$ since $f \in \LtOLtDR,$ i.e.~\cref{con:L} holds.
\epf

\ble[\Cref{con:A} for the Helmholtz stochastic EDP]\label{lem:hh-A}

\noindent If $A \in \LiOLiDRRRdtd,$ $n \in \LiOLiDRRR,$ and $f$ is strongly measurable, then \cref{con:A} holds for the Helmholtz stochastic EDP.
\ele

\bpf[Proof of \cref{lem:hh-A}]
A near-identical argument to that at the beginning of the proof of \cref{lem:hh-L} shows $\toform\circ\coeff$ is strongly measurable. Recall that the Dirichlet-to-Neumann operator $\TrR$ is continuous from $\HhGR$ to $\HmhGR,$ see e.g.~\cite[Theorem 2.6.4]{Ne:01}. Let $\vo \in X, \vt \in Y,$ and observe that the Cauchy--Schwartz inequality and these properties of $\TrR$
imply that there exists $C(k) > 0$ such that 
\begin{align*}
\bigg|\Big[\mleft[\Acomega\mright](\vo)\Big](\vt)\bigg| &= \abs{\int_{\DR} \Big(\mleft(A(\omega) \grad \vo\mright)\cdot\grad \vtb - k^2n(\omega) \vo\vtb \Big)\dd\Leb- \IPGRbig{\DtN \vo}{\vt}}\\
&\hspace{-1cm}\leq \NLiDRspd{A(\omega)}  \NLtDR{\grad \vo}\NLtDR{\grad \vt}\\
&\hspace{-1cm}\quad+ k^2 \NLiDRRR{n(\omega)}\NLtDR{\vo}\NLtDR{\vt}+ C(k) \NHhGR{\gamma\vo}\NHhGR{\gamma\vt},
\end{align*}
where we have used the fact that the two norms
\beq\label{eq:normsdef}
 \NLiDRspd{A(\omega)} = \esssup_{\bx \in \DR} \NopCCd{A(\omega,\bx)} \quad\tand\quad \max_{i,j \in \set{1,\ldots,d}} \NLiDRRR{A_{i,j}(\omega)}
\eeq
are equivalent.
Since the trace operator $\gamma$ is continuous from $\HoDR$ to $\HhGR$ (see, e.g.~\cite[Theorem 3.38]{Mc:00}), there exists $\Ctilde > 0$ such that
\beqs
\Nhomspace{\mleft(\toform\circ\coeff\mright)(\omega)}\leq\Ctilde\max\set{\NLiDRspd{A(\omega)},\NLiDRRR{n(\omega)},C(k)}\Nw{\vo}\Nw{\vt}.
\eeqs
and hence $\toform\circ\coeff \in \LiOhomspace.$
\epf

\subsection{Proofs of \cref{thm:hh-gen,thm:hh-hetero}}\label{sec:applying}

\bpf[Proof of \cref{thm:hh-gen}]
We construct a solution of \cref{prob:msedp} by letting $u = \sol\circ\coeff$ (which is well-defined by \cref{thm:deteu}), and observe that, by construction, $\mleft[a(\omega)\mright]\mleft(u(\omega),v\mright) = \mleft[L(\omega)\mright](v)$ for all $v \in \HozDDR$ almost surely. It follows that $u$ is measurable by \cref{con:hh-fAn,lem:solcont,lem:solcont,lem:contplusmeas}, and so $u$ solves \cref{prob:msedp}. We therefore proceed to apply the general theory.

\Cref{con:coeffstoform,con:coeffstofunc} hold by \cref{lem:hh-AL};
\cref{con:A} holds by \cref{lem:hh-A};
\cref{con:L} holds by \cref{lem:hh-L};
 \cref{con:cborel,con:C} hold by \cref{lem:hh-borelC,con:hh-fAn};
and \cref{con:K} holds by \cref{lem:hh-K}. Therefore we can apply \cref{thm:11,thm:12,lem:svarwelldefined,lem:uniq} to conclude the results.
\epf

\bpf[Proof of \cref{thm:hh-hetero}]
All the conclusions of \cref{thm:hh-gen} hold, and we only need to show that if $u$ solves \cref{prob:msedp} then it also solves \cref{prob:somsedp}. \Cref{con:B} holds by \cref{con:hh-fAn,con:hh-hetero,lem:hh-B}. The result then follows from \cref{thm:3}.
\epf

\appendix

\section{Failure of Fredholm theory for the stochastic variational formulation of Helmholtz problems}\label{sec:federico}
The standard approach to proving existence and uniqueness of a (deterministic) Helmholtz BVP is to  show that the associated sesquilinear form satisfies a G\r{a}rding inequality, and then apply Fredholm theory to deduce that existence and uniqueness are equivalent; see, e.g., \cite[Theorem 4.10]{Mc:00}. This procedure relies on the fact that the inclusion $\HozDDR \hookrightarrow \LtDR$ is compact; see, e.g., \cite[Theorem 3.27]{Mc:00}.

As noted in \cref{sec:otherwork}, the analysis in \cite{FeLiLo:15} of \cref{prob:svsedp} for the Helmholtz Interior Impedance Problem mimics this approach and assums that $\LtOHoD$ is compactly contained  in $\LtOLtD,$ where $D$ is the spatial domain. Here we briefly show $\LtOHoD$ is \emph{not} compactly contained  in $\LtOLtD$ by giving an explicit example of a bounded sequence in  $\LtOHoD$ that has no convergent subsequence in $\LtOLtD.$ Necessary and sufficient conditions for a subset of $L^p\mleft(\mleft[0,T\mright];B\mright),$ for $B$ a Banach space, to be compact, can be found in \cite{Si:86}. In particular, \cite{Si:86} shows that a space $C$ being compactly contained in a space $B$ does not by itself imply $L^2\mleft(\mleft[0,T\mright];C\mright)$ is compactly contained in $L^2\mleft(\mleft[0,T\mright];B\mright).$

\begin{example} \label{ex:federico}
Let $(\Omega,\cF,\PP) = (\mleft[0,1\mright],\cB(\mleft[0,1\mright]),\lambda).$ Let $D$ be a compact subset of $\mathbb{R}^d.$ Since $\LtO$ is separable, it has an orthonormal basis, which we denote by $(f_m)_{m \in \NN}.$ Let $u_m \in  \LtOHoD$ be defined by $ u_m(\omega)(x) \de f_m(\omega), \tfa x \in D,$
i.e., for each value of $\omega,$ $u_m(\omega)$ is a constant function on $D$ and so $\NHoD{u_m(\omega)} = \NLtD{u_m(\omega)}.$ Then
\beqs
\NLtOHoD{u_m}^2 = \int_\Omega \NHoD{u_m(\omega)}^2 \dd\PP(\omega) = \lambda(D)^2\int_\Omega \abs{f_m(\omega)}^2 \dd\PP(\omega)= \NLtO{f_m}^2 \lambda(D)^2,
\eeqs
and so $u_m$ is a bounded sequence in $\LtOHoD.$ However, for $n \neq m,$ we have
\begin{align*}
\NLtOLtD{u_m-u_n}^2 &= \int_\Omega \NLtD{u_m(\omega)-u_n(\omega)}^2 \dd\PP(\omega)\\
&= \lambda(D)^2 \int_\Omega \abs{u_m(\omega) - u_n(\omega)}^2 \dd\PP(\omega) = \lambda(D)^2\NLtO{f_m-f_n}^2 = 2\lambda(D)^2
\end{align*}
if $n \neq m,$ since the $f_m$ form an orthonormal basis for $\LtD.$ Therefore $(u_m)_{m \in \NN}$ is bounded in $\LtOHoD$ but does not have a convergent subsequence in $\LtOLtD,$ and thus the inclusion of $\LtOHoD$ into $\LtOLtD$ cannot be compact.
\end{example}

\section{Recap of basic material on measure theory and Bochner spaces}\label{app:mtBs}
We include this section, not only for completeness, but also to aid readers of this paper who are more familiar with deterministic, as opposed to stochastic, Helmholtz problems. Recall that here, and in the rest of the paper, $\OFP$ is a complete probability space.

\subsection{Recap of measure theory results}

We first recall some results from measure theory, with our main reference \cite{Bo:07}. Even though \cite{Bo:07} mainly considers maps with image $\RR,$ the results we quote for more general images are straightforward generalisations of the results in \cite{Bo:07}.

\begin{definition}[Measurable map]\label{def:meas}
If $(\M,\FM)$ and $(N,\FN)$ are measurable spaces, we say that $f:\M\rightarrow N$ is measurable (with respect to $(\FM,\FN)$) if $f^{-1}(E) \in \FM$ for all $E \in \FN.$
\end{definition}

\bde[Borel {$\sigma$}-algebra]\label{def:borelsigma}
If $(S,\TS)$ is a topological space, the \defn{Borel $\sigma$-algebra} $\Borel{S}$ on $S$ is the $\sigma$-algebra generated by $\TS.$
\ede

If $V$ is any topological space (including a Hilbert, Banach, metric, or normed vector space) then we will take always the Borel $\sigma$-algebra on $V$ unless stated otherwise.

\ble[Continuous maps are measurable { \cite[Theorem 2.1.2]{Bo:07}}] \label{lem:contmeas}
Any continuous function between two topological spaces is measurable.
\ele

\ble[The composition of a measurable and a continuous map is measurable {\cite[Text at top of p. 146]{Bo:07}}]\label{lem:contplusmeas}
Let $(\M,\FM)$ be a measurable space and let $\mleft(S,\TS\mright)$ and $\mleft(T,\TopT\mright)$ be topological spaces. Let $f:M \rightarrow S$ be measurable and let $h : S \rightarrow T$ be continuous. Then $h \circ f$ is measurable.
\ele

\bde[Product $\sigma$-algebra {\cite[Section IV.11]{Do:94}}]\label{def:prodsigma}
Let $\mleft(\M_1,\FM_1\mright),\ldots,\mleft(\M_m,\FM_m\mright)$ be measurable spaces. The \emph{product $\sigma$-algebra} $\M_1\otimes\cdots\otimes\M_m$ is defined as the $\sigma$-algebra generated by the set of measurable rectangles
$
\set{R_{1} \times \cdots \times R_m  \st R_{1} \in \FM_1, \ldots, R_m \in \FM_m}.
$
\ede

\ble[Measurability of the Cartesian product of measurable functions]\label{lem:measprod}

Let $\mleft(\M_1,\FM_1\mright),\ldots,\mleft(\M_m,\FM_m\mright)$ be measurable spaces and $h_j:\Omega \rightarrow \M_j,\, j = 1,\ldots,m$ be measurable functions. Then the product map $\Prodf:\Omega \rightarrow \M_1 \times \cdots \times \M_m$ given by
$\Prodf(\omega) \de \mleft(h_1(\omega),\ldots,h_m(\omega)\mright)$
is measurable with respect to $\mleft(\cF,\FM_1 \otimes \cdots\otimes\FM_m\mright).$
\ele

\bpf[Sketch proof of \cref{lem:measprod}]

Let $\measrectmany{\FM_1,\ldots,\FM_m}$ denote the set of measurable rectangles, as in \cref{def:prodsigma}. Let
$\cP \de \set{C \subseteq \M_1 \times \cdots \times \M_m \st \Prodf^{-1}\mleft(C\mright) \in \cF}.$
The proof of the lemma consists of the following straightforward steps, whose proofs are omitted:
(i) Show 

\noindent $\measrectmany{\FM_1,\ldots,\FM_m} \subseteq \cP.$
(ii) Show $\cP$ is a $\sigma$-algebra.
(iii) Deduce $\FM_1 \otimes \cdots \otimes \FM_m \subseteq \cP$ (since $\FM_1 \otimes \cdots \otimes \FM_m$ is generated by measurable rectangles).
(iv) Conclude $\Prodf$ is measurable with respect to $\mleft(\cF,\FM_1\otimes\cdots\otimes\FM_m\mright).$
\epf

\ble[Product of Borel $\sigma$-algebras is Borel $\sigma$-algebra of the product {\cite[Lemma 6.2.1 (i)]{Bo:07}}]\label{lem:bogachev}
Let $H_1,H_2$ be Hausdorff spaces and let $H_2$ have a countable base (e.g.~$H_2$ could be a separable metric space). Then $\Borel{H_1\times H_2} = \Borel{H_1} \otimes \Borel{H_2},$ where $\Borel{H_1\times H_2}$ is the Borel $\sigma$-algebra of the product topology on $H_1\times H_2.$
\ele

\subsection{Recap of results on Bochner spaces}

We now recap the theory of Bochner spaces, using \cite{DiUh:77} as our main reference. In what follows the space $V$ is always a Banach space.

\bde[Simple function]
A function $v:\Omega \rightarrow V$ is \defn{simple} if there exist $v_1,\ldots,v_m \in V$ and $E_1,\ldots,E_m \in \cF$ such that
$v = \sum_{i=1}^m v_i \chi_{E_{i}},$
where $\chi_{E_{i}}$ is the indicator function on $E_{i}.$
\ede

\bde[Strongly measurable]\label{def:strongmeas}
A function $v:\Omega \rightarrow V$ is \defn{strongly measurable}
\footnote{In \cite{DiUh:77} the authors use the term \defn{$\mu$-measurable} instead of \defn{strongly measurable} (where $\mu$ is the measure on the domain of the functions under consideration).} 
if there exists a sequence of simple functions $(\vn)_{n \in \NN}$ such that
$\lim_{n \rightarrow \infty} \NV{\vn - v} = 0$, 
$\PP$-almost everywhere.
\ede

\bde[Bochner integrable {\cite[p. 49]{DiUh:77}}]
A strongly measurable function $v:\Omega \rightarrow V$ is called \defn{Bochner integrable} if there exists a sequence of simple functions $(\vn)_{n \in \NN}$ such that
$\lim_{n \rightarrow \infty} \int_{\Omega} \NV{\vn(\omega) - v(\omega)} \dd\PP(\omega) = 0.$
\ede

\bth[Condition for Bochner integrability {\cite[Theorem II.2.2]{DiUh:77}}]\label{thm:bochnercond}
A strongly measurable function $v:\Omega \rightarrow V$ is Bochner integrable if and only if $\int_\Omega \NV{v} \dd\PP < \infty.$
\enth

\bco[Sufficient condition for Bochner integrability]\label{cor:bochnersimple}
Let $p \geq 1.$ If a strongly measurable function $v:\Omega \rightarrow V$ has $\int_\Omega \NV{v}^p \dd\PP < \infty,$ then $v$ is Bochner integrable.
\eco

\bde[Bochner norm]\label{def:bochnernorm}
For a Bochner integrable function $v:\Omega\rightarrow V,$ let
\beqs
\NLpOV{v} \de \mleft(\int_\Omega \NV{v(\omega)}^p \dd\PP(\omega)\mright)^{1/p}, \,1 \leq p < \infty, \,\,\text{and}\,\,\NLiOV{v} \de \esssup_{\omega \in \Omega} \NV{v(\omega)}.
\eeqs
\ede

\bde[Bochner space]\label{def:bochnerspace}
Let $1\leq p \leq \infty.$ Then
\beqs
\LpOV \de \set{v:\Omega\rightarrow V \st v \text{ is Bochner integrable,}\,\NLpOV{v} < \infty}.
\eeqs
\ede

\bde[Complete probability space]
A probability space $\OFP$ is complete if for every $\Eo \in \cF$ with $\PP(\Eo)=0,$ the inclusion $\Et \subseteq \Eo$ implies that $\Et \in \cF.$
\ede

\bde[Separable space]
A topological space is \defn{separable} if it contains a countable, dense subset.
\ede

\bde[$\sigma$-finite]
A probability space $\OFP$ is \defn{$\sigma$-finite} if there exist $E_1,E_2,\ldots \in \cF$ with $\PP(E_m) < \infty$ for all $m \in \NN$ such that $\Omega = \cup_{m=1}^\infty E_m.$
\ede

\bth[Pettis measurability theorem {\cite[Proposition 2.15]{Ry:02}}]\label{thm:pettis}
Let $\OFP$ be a complete $\sigma$-finite measure space. The following are equivalent for a function $v:\Omega \rightarrow V$:
(i) $v$ is  strongly measurable,
(ii) $v$ is measurable and $\PP$-essentially separably valued.
\enth

\bco[Equivalence of measurable and strongly measurable when the image is separable]\label{cor:pettis}
Let $\OFP$ be a complete $\sigma$-finite measure space. If $V$ is a separable Banach space, then a function $v:\Omega\rightarrow V$ is strongly measurable if, and only if, it is measurable.
\eco

\ble[The composition of a continuous map and a {$\PP$}-essentially separably valued map]\label{lem:esssep}
Let $\mleft(S,\TS\mright)$ and $\mleft(T,\TopT\mright)$ be topological spaces. If $f_1:\Omega \rightarrow S$ and $f_2:S\rightarrow T$ are such that $f_1$ is $\PP$-essentially separably valued and $f_2$ is continuous, then $f_2\circ f_1$ is $\PP$-essentially separably valued.
\ele

\bpf[Proof of \cref{lem:esssep}]
As $f_1$ is $\PP$-essentially separably valued, there exists $E \in \cF$ such that $\PP(E) = 1$ and $f_1(E) \subseteq G \subseteq S,$ where $G$ is separable. As $\ft$ is continuous, $\ft(G)$ is separable \cite[Theorem 16.4(a)]{Wi:70}. Therefore, since $\mleft(\ft \circ\fo\mright)(E) \subseteq \ft(G),$ it follows that $\ft\circ\fo$ is $\PP$-essentially separably valued.
\epf

\ble[The composition of a continuous map and a strongly measurable map]\label{lem:contplusstrong}
If $\Bo$ and $\Bt$ are Banach spaces and there exist $f_1:\Omega \rightarrow \Bo$ and $f_2:\Bo\rightarrow \Bt$ such that $f_1$ is strongly measurable and $f_2$ is continuous, then $f_2\circ f_1$ is strongly measurable.
\ele

\bpf[Proof of \cref{lem:contplusstrong}]
By \cref{thm:pettis}, $\fo$ is both measurable and $\PP$-essentially separably valued. Therefore we can apply \cref{lem:contplusmeas,lem:esssep} to conclude $\ft \circ \fo$ is both measurable and $\PP$-essentially separably valued. Hence by \cref{thm:pettis} $\ft\circ\fo$ is strongly measurable.
\epf

\ble[Zero in all integrals implies zero almost everywhere {\cite[Corollary II.2.5]{DiUh:77}}]\label{lem:gotoae}
If $\alpha$ is  Bochner integrable and  $\int_E \alpha(\omega) \dd\PP(\omega) = 0 $ for each $E \in \cF$ then $\alpha=0$ $\PP$-almost everywhere.
\ele

\ble[Cartesian product of $\PP$-essentially separably valued maps]\label{lem:prodsep}
Let 

\noindent $\mleft(\cCo,\Top{\cCo}\mright),\ldots,\mleft(\cCm,\Top{\cCm}\mright)$ be topological spaces, and let $s_j:\Omega\rightarrow\cCj,\,j=1,\ldots,m$ be $\PP$-essentially separably valued. Define $\cC \de \cCo \times \cdots \times \cCm$ and equip $\cC$ with the product topology. Then the map $f:\Omega\rightarrow \cC$ given by $s(\omega) \de \mleft(s_1(\omega),\ldots,s_m(\omega)\mright) $ is $\PP$-essentially separably valued.
\ele

The proof of \cref{lem:prodsep} is straightforward and omitted.

\section{Measurability of series expansions (used in \cref{sec:generating})}\label{app:meas}

Here we collect together results from measure theory that allow us to conclude in \cref{lem:seriesmeas} that the series expansions for $A$ and $n$ in \cref{sec:generating} are measurable. As mentioned in \cref{sec:generating}, the proof that the sum of measurable functions is measurable is standard, but we have not been able to find this result stated in the literature for this particular setting of mappings into a separable subspace of a general normed vector space.
\ble\label{lem:sepsum}
If $U$ is a separable normed vector space, $m \in \NN,$ and $\phi_j:\Omega\rightarrow U,$ $j=1,\ldots,m$  are measurable functions, then $\phi_1+\cdots+\phi_m : \Omega\rightarrow U$ is measurable.
\ele

\bpf[Sketch proof of \cref{lem:sepsum}]
By induction, it is sufficient to show the result for $m=2.$ We let $\Ball{U}{r}{v}$ denote the ball of radius $r>0$ about $v \in U$. To show $\phio+\phit$ is measurable, we let $v \in U, r>0$ and we show $\mleft(\phio+\phit\mright)^{-1}\mleft(\Ball{U}{r}{v}\mright) \in \cF.$ Let $\QQU$ denote a countable dense subset of $U,$ which exists as $U$ is separable. Let $\QQFF$ denote a countable dense subset of the field $\FF$, which exists as $\FF = \RR$ or $\CC.$

For $s \in \QQFF, q \in \QQU$ let
\beqs
\setsq = \set{\omega \in \Omega \st \NU{\phio(\omega)-\half v - q} < s} \cap \set{\omega \in \Omega \st \NU{\phit(\omega)-\half v + q} < r-s}.
\eeqs
We claim
\beq\label{eq:sumseteq}
\mleft(\phio+\phit\mright)^{-1}\mleft(\Ball{U}{r}{v}\mright) = \bigcup_{s \in \QQFF} \bigcup_{q \in \QQU} \setsq,
\eeq
and the result then follows as the right-hand side is an element of the $\sigma$-algebra $\cF.$ To show \eqref{eq:sumseteq}, let $\omega \in \cup_{s \in \QQFF}\cup_{q \in \QQU} \setsq,$ and let $s \in \QQFF, q \in \QQU$ be such that $\omega \in \setsq.$ Then it follows from the triangle inequality that $\omega \in \mleft(\phio+\phit\mright)^{-1}\mleft(\Ball{U}{r}{v}\mright).$
Now let $\omega \in \mleft(\phio+\phit\mright)^{-1}\mleft(\Ball{U}{r}{v}\mright),$ define $\dw \de r - \NU{\phio(\omega)+\phit(\omega)-v} > 0,$ fix $s \in \QQFF \cap (0,\dw/2),$ and choose $q \in \QQU$ such that $\NU{\phio(\omega) - v/2 -q} < s.$ Then again it follows from the triangle inequality that $\omega \in \setsq,$ and thus \eqref{eq:sumseteq} holds, as required.
\epf

\bco\label{cor:sepsubsum}
If $V$ is a normed vector space, $U \subseteq V$ is a separable subspace, and $\phi_j:\Omega\rightarrow U,$  $j=1,\ldots,m$  are measurable functions, then $\phi_1+\cdots+\phi_m : \Omega\rightarrow U$ is measurable.
\eco

\ble\label{lem:scalarmultmeas}
Let $V$ be a normed vector space. If $v \in V$ and $Y:\Omega\rightarrow \FF$ is a measurable function, then $Yv:\Omega\rightarrow V$ is a measurable function.
\ele

\bpf[Proof of \cref{lem:scalarmultmeas}]
The map $\scalarmult:\FF\rightarrow V$ given by $\scalarmult(x) = xv$ is continuous. As $Yv = \scalarmult \circ Y,$ it follows from \cref{lem:contplusmeas} that $Yv$ is measurable.
\epf

\ble\label{lem:incborel}
If $V$ is a normed vector space and $U \subseteq V,$ then the inclusion map $\inc:U\rightarrow V$ is measurable.
\ele

\bpf[Proof of \cref{lem:incborel}]
As $\inc$ is continuous, it immediately follows that it is measurable.
\epf

\bco\label{cor:meassubmeansmeas}
If $V$ is a normed vector space, $U\subseteq V$ and $\phi:\Omega\rightarrow U$ is measurable, then $\phi:\Omega\rightarrow V$ is measurable.
\eco

\bpf[Proof of \cref{cor:meassubmeansmeas}]
This is immediate from \cref{lem:incborel} and \cref{lem:contplusmeas}.
\epf

\ble\label{lem:spansep}
If $V$ is a normed vector space, $m \in \NN,$ and $\phi_1,\ldots,\phi_m \in V$ for $j = 1,\ldots,m$ then $\spanset{\phi_1,\ldots,\phi_m}$ is a separable subspace of $V.$
\ele

\bpf[Sketch Proof of \cref{lem:spansep}]
As $\FF = \RR$ or $\CC,$ it has a separable subset $\QQFF.$ Since a finite product of countable sets is countable, the set
\beqs
\set{\Ball{V}{1/n}{q_1\phi_1 + \cdots + q_m \phi_m} \st n \in \NN, q_1,\ldots,q_m \in \QQFF}
\eeqs
is a countable base for the topology on $\spanset{\phi_1,\ldots,\phi_m}$ induced by the norm $\NV{\cdot}.$ 
\epf

\ble\label{lem:summeas}
If $V$ is a normed vector space, $m \in \NN,$ and for $j = 1,\ldots,m$, $\phi_j\in V$ and $Y_j : \Omega \rightarrow \FF$ are measurable, then the function $\phi:\Omega\rightarrow V$ given by
\beqs
\phi(\omega) = \phi_0 + \sum_{j=1}^m Y_j(\omega)\phi_j
\eeqs
is measurable.
\ele

\bpf[Proof of \cref{lem:summeas}]
The subspace $U=\spanset{\phi_0,\phi_1,\ldots,\phi_m}$ is separable by \cref{lem:spansep}, and it is clear that the image of $\phi$ lies in $U.$ By \cref{lem:scalarmultmeas} and \cref{cor:sepsubsum}, $\phi:\Omega\rightarrow U$ is measurable, and therefore $\phi:\Omega\rightarrow V$ is measurable by \cref{cor:meassubmeansmeas}.
\epf

We now prove that almost-surely convergent sequences of measurable functions are measurable, and we then apply this result to the partial sums in the definitions of $A$ and $n$ in \eqref{eq:nseries}.

We will use the following theorem to establish that the almost-sure limit of a sequence of measurable functions is measurable.

\bth[{\cite[Theorem 4.2.2]{Du:02}}]\label{thm:dudley}
Let $(\metsp,\metric)$ be a metric space. Suppose the functions $\zeta_j:\Omega \rightarrow \metsp$ are measurable, for all $j \in \NN.$ If the limit
\beqs
\zeta(\omega) = \lim_{j \rightarrow \infty} \zeta_j(\omega)
\eeqs
exists for every $\omega \in \Omega,$ then the function $\zeta:\Omega \rightarrow \metsp$ is measurable.
\enth

\bco\label{lem:paelimitmeas}
Let $(\metsp,\metric)$ be a metric space. Suppose the functions $\zeta_m:\Omega \rightarrow \metsp$ are measurable, for all $m \in \NN.$ If the limit
\beq\label{eq:zetalimit}
\lim_{m \rightarrow \infty} \zeta_m(\omega)
\eeq
exists almost surely, then there exists a measurable function $\zeta:\Omega \rightarrow \metsp$ such that
\beqs
\zeta(\omega) = \lim_{m \rightarrow \infty} \zeta_m(\omega)
\eeqs
whenever the limit exists.
\eco

\bpf[Proof of \cref{lem:paelimitmeas}]
Following \cite{El:11}, we define $\Omegatilde = \set{\omega \in \Omega \st \text{ \eqref{eq:zetalimit} exists}}.$ Then, for $m \in \NN$ define $\zetatilde_m:\Omega\rightarrow \metsp$ by
\beqs
\zetatilde_m(\omega) =
\begin{cases}
\zeta_m(\omega) & \tif \omega \in \Omegatilde\\
0 & \tif \omega \not\in\Omegatilde
\end{cases}
\eeqs
Observe that, by construction, the limit $\zetatilde(\omega) = \lim_{m\rightarrow\infty} \zetatilde_m(\omega)$ exists \emph{for all} $\omega \in \Omega$ and the functions $\zetatilde_m$ are measurable. Therefore, by \cref{thm:dudley}, $\zetatilde$ is measurable.
\epf

\ble\label{lem:paeseriesmeas}
Let $V$ be a normed vector space. If there exist $\phi_j\in V,$ $j=0,1,\ldots$ and measurable functions $\Yj:\Omega \rightarrow \FF,$ $j \in \NN$ such that the series
\beqs
\phi_0 + \sum_{j=1}^\infty \Yj(\omega)\phi_j
\eeqs
exists in $V$ almost surely, then there exists a measurable function $\phi:\Omega\rightarrow V$ such that
\beqs
\phi(\omega) = \phi_0 + \sum_{j=1}^\infty \Yj(\omega)\phi_j
\eeqs
almost surely.
\ele

\bpf[Proof of \cref{lem:paeseriesmeas}]
By \cref{lem:summeas}, the partial sums $\phi_0 + \sum_{j=1}^m \Yj(\omega)\phi_j$, for $m \in \NN$ are measurable, and by assumption their limit as $m \rightarrow \infty$ exists almost surely. Therefore, applying \cref{lem:paelimitmeas} to the partial sums, we obtain the result.
\epf

\ble\label{lem:seriesexistsas}
The series expansions for both $A$ and $n$ defined by \eqref{eq:nseries} exist in $\WoiDRRRdtd$ and $\WoiDRRR$ almost surely, respectively.
\ele

\bpf[Proof of \cref{lem:seriesexistsas}]
The spaces $\WoiDRRRdtd$ and $\WoiDRRR$ are Banach spaces, by definition of their norms (see, e.g., \eqref{eq:normsdef}). Therefore it suffices to show that the partial sums of the series expansions for $A$ and $n$ in \eqref{eq:nseries} are Cauchy sequences. As the proofs for $A$ and $n$ are completely analogous, we only give the proof for $A$ here.

First observe that since each of the random variables $\Yj$ in \eqref{eq:nseries} is uniformly distributed on $[-1/2,1/2],$ it follows that for all $j \in \NN$, $\esssup_{\omega \in \Omega} \abs{\Yj(\omega)} = \frac12.$
Therefore, we can conclude that the bound $\esssup_{\omega \in \Omega} \sup_{j \in \NN} \abs{\Yj(\omega)} \leq \half$ holds.
(For if not, then, there would exist $\Omegahat \subseteq \Omega$ with $\PP(\Omegahat) > 0$ such that for all $\omega \in \Omegahat$, $\sup_{j \in \NN} \abs{\Yj(\omega)} > \half.$
Then there would exist $\jhat \in \NN$ such that $\abs{\Yjhat(\omega)} > 1/2$ for all $\omega \in \Omegahat,$ which would give the contradiction $\esssup_{\omega \in \Omega} |\Yjhat(\omega)| > \half$.)

It now suffices to show that for $\PP$-almost every $\omega \in \Omega,$ the partial sums of the series expansion in \eqref{eq:nseries} form a Cauchy sequence. Recall that for $\PP$-almost every $\omega \in \Omega$
\beqs
\sup_{j \in \NN} \abs{\Yj(\omega)} \leq \half.
\eeqs
For such an $\omega$, and $m \in \NN,$ define the $m$th partial sum
\beqs
\Am(\omega) = \Az + \sum_{j=1}^m \Yj(\omega)\Psij.
\eeqs
It is straightforward to show that $\mleft(\Am(\omega)\mright)_{m \in \NN}$ is a Cauchy sequence in $\WoiDRRRdtd$, using the assumption \eqref{eq:Apsimeas}; therefore, the series expansion for $A(\omega)$ in \eqref{eq:nseries} exists almost surely.
\epf

\ble\label{lem:seriesmeas}
The functions $A$ and $n$ defined by \eqref{eq:nseries} are measurable.
\ele

\bpf[Proof of \cref{lem:seriesmeas}]
The result immediately follows from \cref{lem:seriesexistsas,lem:paeseriesmeas}.
\epf

\section*{Acknowledgements}
We thank Federico Cornalba (University of Bath) for the example in \cref{sec:federico} and for introducing us to \cite{Si:86}. We also thank Jack Betteridge (Bath), Ivan Graham (Bath),  Kieran Jarrett (Bath), Robert Scheichl (Bath/Universit\"at Heidelburg), and Tony Shardlow (Bath) for useful discussions. We thank the anonymous referees for their comments, which improved the organisation of the paper.


\bibliographystyle{siamplain}
\bibliography{biblio_svar}

\end{document}